\documentclass{amsart}
\usepackage{amssymb}
\usepackage{amsfonts}
\usepackage{latexsym}

\newtheorem{theorem}{Theorem}[section]
\newtheorem{lemma}[theorem]{Lemma}
\newtheorem{proposition}[theorem]{Proposition}
\newtheorem{corollary}[theorem]{Corollary} 
\theoremstyle{definition}  
\newtheorem{definition}[theorem]{Definition}

\newtheorem{example}[theorem]{Example}

\newtheorem{conjecture}[theorem]{Conjecture}  
\newtheorem{remark}[theorem]{Remark}

\newcommand{\Tr}{\text{Tr}}

\newcommand{\Ker}{\text{Ker\,}}

\newcommand{\End}{\text{End}} 

\newcommand{\Hom}{\text{Hom}}

\newcommand{\Ind}{\text{Ind}}
\newcommand{\Rep}{\text{Rep}}

\newcommand{\eps}{\varepsilon}

\newcommand{\ben}{\begin{enumerate}}
\newcommand{\een}{\end{enumerate}}

\newcommand{\Cl}{{\text{Cl}}}
\newcommand{\mC}{{\mathcal C}}

\newcommand{\mD}{{\mathcal D}}

\newcommand{\mM}{{\mathcal M}}

\newcommand{\cC}{{\mathcal C}}
\newcommand{\cA}{{\mathcal A}}

\newcommand{\fA}{{\mathfrak A}}
\newcommand{\fB}{{\mathfrak B}}
\newcommand{\cD}{{\mathcal D}}
\newcommand{\cM}{{\mathcal M}}
\newcommand{\be}{{\bf 1}}
\newcommand{\p}{\partial}

\newcommand{\BZ}{{\mathbb Z}}

\newcommand{\Ve}{\mbox{Vec}}
\newcommand{\chara}{\mbox{char}}
\newcommand{\Fun}{\mbox{Fun}}
\newcommand{\Mod}{\mbox{Mod}}
\newcommand{\Bimod}{\mbox{Bimod}}
\newcommand{\iHom}{\underline{\mbox{Hom}}}

\begin{document}

\title{Finite tensor categories}
\begin{abstract} We start the general structure theory 
of not necessarily semisimple finite tensor categories, 
generalizing the results in the semisimple case (i.e. for fusion
categories), obtained recently in our joint work with
D. Nikshych. In particular, we generalize 
to the categorical setting the Hopf and quasi-Hopf algebra freeness
theorems due to Nichols-Zoeller and Schauenburg, respectively.
We also give categorical versions of the theory of 
distinguished group-like elements in a finite dimensional Hopf
algebra, of Lorenz's result on degeneracy of the Cartan matrix, 
and of the absence of primitive elements in a finite dimensional Hopf
algebra in zero characteristic. We also develop the theory of
module categories and dual categories for not necessarily semisimple
finite tensor categories; the crucial new notion here is that of an 
{\it exact module category}. Finally, we classify indecomposable
exact module categories over the simplest finite tensor
categories, such as representations of a finite group in positive
characteristic, representations of a finite supergroup, and 
representations of the Taft Hopf algebra. 
\end{abstract}

\author{Pavel Etingof}
\address{Department of Mathematics, Massachusetts Institute of Technology,
Cambridge, MA 02139, USA}
\email{etingof@math.mit.edu}

\author{Viktor Ostrik}
\address{Department of Mathematics, Massachusetts Institute of Technology,
Cambridge, MA 02139, USA}
\email{ostrik@math.mit.edu}

\maketitle

{\bf Dedicated to Boris Feigin on the occasion of his 50-th birthday}
\vskip .1in

\section{Introduction}

The aim of this paper is to develop a systematic 
theory of not necessarily semisimple finite tensor and
multi-tensor categories, 
similarly to how it was done in \cite{ENO} and references therein
in the semisimple case, i.e., for fusion and multi-fusion categories. 
There are several (interrelated) motivations for this: 

1. Representations of finite groups in positive characteristic. 

2. Finite dimensional Hopf algebras, in particular quantum groups 
$U_q({\frak g})$ at roots of unity.

3. Logarithmic conformal field theories; they lead to 
nonsemisimple finite tensor categories, similarly to how 
rational conformal field theories lead to semisimple ones (see \cite{Ga}). 

4. Fusion categories of zero global dimension (duals to such categories 
may be nonsemsimple).

We begin by studying the general properties of finite tensor categories,
focusing on issues specific to the nonsemisimple situation,
like the behavior of projective objects (see Section 2). 
More specifically, we generalize a number of classical results in the theory 
of finite dimensional Hopf algebras to the categorical setting. 
For example, we show that a surjective quasi-tensor 
(in particular, tensor) functor\linebreak  
$\mC\to \mD$ between 
finite tensor categories maps projective objects to projective ones,
and that the regular (virtual) object of $\mC$ maps to a multiple
of the regular object of $\mD$. This means that the 
Frobenius-Perron dimension of $\mC$ is divisible by that of $\mD$,
and implies as a special case the Hopf and quasi-Hopf algebra freeness theorems
of Nichols-Zoeller and Schauenburg, respectively. We also generalize
to the categorical setting the theory of distinguished grouplike 
elements for finite dimensional Hopf algebras, Lorenz's theorem on the 
degeneracy of the Cartan matrix, and the theorem that a finite dimensional
Hopf algebra in characteristic zero cannot have nonzero primitive elements. 
This last generalization implies that any finite tensor category 
in zero characteristic with a unique simple object is equivalent 
to the category of vector spaces. This generalizes the fact 
that a local finite dimensional Hopf algebra in zero characteristic 
is 1-dimensional. 

More significantly, in Section 3 we propose a generalization to the 
nonsemisimple case
of the theory of module categories and dual categories (\cite{O},\cite{ENO}). 
The naive generalization does not give a good theory: 
for example, if $\mC$ is the category of vector spaces and $A$ any finite 
dimensional algebra, then $\mM=\Rep(A)$ is a module category 
over $\mC$, so there is no hope of explicit classification of 
module categories even over the simplest possible tensor category
(the category of vector spaces). The situation with dual categories
is even worse: the dual category $\mC_{\mM}^*$ (i.e. the category 
of $\mC$-linear functors from $\mM$ to itself) is the category 
of $A$-bimodules with the bimodule tensor product. 
If the algebra $A$ is not semisimple, this category is not rigid
and the tensor product functor in it is not exact, so much of the theory 
fails. 
This shows tha one should not study all module categories, 
but rather restrict to a ``correct'' subclass of them, 
containing some desirable examples, such as

1) any finite tensor category $\mC$ as a module over itself;

2) semisimple module categories (= weak Hopf algebras);

3) any finite tensor category $\mC$ as a module over $\mC\boxtimes \mC^{op}$
(in this case the dual is the Drinfeld center $Z(\mC)$ of $\mC$). 
 
In this paper, we propose such a subclass. Namely, we define an {\em exact} 
module category over $\mC$ to be any module category $\mM$ such that 
for any object $X\in \mM$ and any {\it projective} object 
$P\in \mC$ the product $P\otimes X$ is projective. This class 
contains examples 1-3, and reduces to example 2 for semisimple $\mC$. 
We show that for exact module categories 
the theory works as perfectly as it does in the semisimple case 
(for fusion and multifusion categories). In particular, we show that 
if $\mC$ is a finite tensor category and $\mM$ an indecomposable exact 
module category over $\mC$, then $\mC^*_{\mM}$ is a 
finite tensor category of the same Frobenius-Perron dimension, 
and $(\mC_\mM^*)_\mM^*=\mC$. In particular, $Z(\mC)$ is a finite tensor 
category, whose Frobenius-Perron dimension is the square of that of $\mC$. 

Finally, in Section 4 we classify exact module categories 
over three examples of finite tensor categories: 

1) representations of a finite group in positive characteristic,

2) representations of a finite supergroup in characteristic $\ne 2$, and

3) representations of the Taft Hopf algebra. 

In example 1, our result is a generalization to positive characteristic
of the result of \cite{Ob}; in example 2, it is a generalization
of the work \cite{EG1}.  

{\bf Acknowledgements.}
It is our joy to dedicate this paper to the 50-th birthday 
of Boris Feigin. Pasha Etingof is glad to use this opportunity 
to express his deep gratitude to Borya for giving him
mathematical education in 1986-1990, which enabled him to do mathematics.

We thank J. de Jong, S. Gelaki, D. Nikshych, and J. Starr 
for useful discussions and references. 
The research of P.E. was partially supported by
the NSF grant DMS-9988796, and was done in part for the Clay Mathematics 
Institute. The research of V.O. was supported by the NSF 
grant DMS-0098830.  

\section{General properties of finite tensor categories}

\subsection{Definitions and notation}
Let $k$ be an algebraically closed field. Let $\mC$ be an abelain
category over $k$, where morphism spaces are finite dimensional, and
every object has finite length. We will say that $\mC$ is {\it finite} if it
has finitely many simple objects, and each of them has a projective
cover (we will denote the projective cover of a simple object $X\in \cC$ by
$P(X)$). This is equivalent to $\mC$ being equivalent to the category of
finite dimensional representations of a finite dimensional $k$-algebra. 

By a {\it tensor} category we will mean an abelian rigid
tensor category over $k$ in which the unit object $\bold 1$ is
simple (see \cite{BaKi} for a full definition). 
It is known (\cite{BaKi}, Proposition 2.1.8) that in such a
category, the tensor product functor is exact 
in both arguments. 

The main object of study in this paper is {\it finite tensor categories.}
For example, if $H$ is a finite dimensional Hopf (or, more generally, 
quasi-Hopf)
algebra over $k$, then $\Rep H$ is a finite tensor category. 

Let $\mC$ be a finite tensor category, $I$ be the set 
of isomorphism classes of simple objects of $\mC$, and let $i^*,{}^*i$ denote 
the right and left duals to $i$, respectively. 
Let $Gr(\mC)$ be the Grothendieck ring of $\mC$, 
spanned by isomorphism classes of the simple objects $L_i$. 
In this ring, we have $L_iL_j=\sum_k{N_{ij}^k}L_k$, where $N_{ij}^k$ are 
positive integers. Also, let $P_i$ be the projective covers of $L_i$. 

We will use the symbol $\boxtimes$ for Deligne's tensor product of abelian
categories, see \cite{De2}. Recall that for two finite dimensional algebras
$A$ and $B$ one has $\Rep(A)\boxtimes \Rep(B)=\Rep(A\otimes B)$, see 
\cite{De2}. Note that if $\cC$ and $\cD$ are finite tensor
categories then $\cC \boxtimes \cD$ also has a natural structure of a 
finite tensor category. 

\subsection{Projectivity of tensor products and duals} The following
Proposition is well known, see e.g. \cite{KL4}, p. 441, Corollary 2.

\begin{proposition}\label{tensprod}
Let $P$ be a projective object in $\mC$, and $X$ any object of $\mC$. Then 
$P\otimes X$ is projective. 
\end{proposition}

\begin{proof}
We have $\Hom(P\otimes X,Y)=\Hom(P,Y\otimes X^*)$. 
The latter functor in $Y$ is exact, since the tensor product is biexact. 
Thus, $P\otimes X$ is projective. 
\end{proof} 

Similarly, $X\otimes P$ is projective for a projective $P$. 

Let $[Z:L]$ be the multiplicity of occurence of a simple object $L$ in 
an object $Z$. 

\begin{proposition}\label{tensprod1} 
For any object $Z$ of $\mC$, 
$P_i\otimes Z\cong \oplus_{j,k} N_{kj^*}^i [Z:L_j]P_k$, and 
$Z\otimes P_i\cong \oplus_{j,k} N_{{}^*jk}^i [Z:L_j]P_k$.
\end{proposition}

\begin{proof}
$\Hom(P_i\otimes Z,L_k)=\Hom(P_i,L_k\otimes Z^*)$, and the first formula 
follows from Proposition \ref{tensprod}. 
The second formula is analogous. 
\end{proof}

\begin{proposition} \label{dualproj} Let $P$ be a projective object in $\mC$. 
Then $P^*$ is also projective. 
\end{proposition}

\begin{proof}
We need to show that the functor $\Hom(P^*,?)$ is exact. This functor
is isomorphic to $\Hom(\be,P\otimes ?)$. The functor $P\otimes ?$ is exact
and moreover any exact sequence splits after tensoring with $P$ by
Proposition \ref{tensprod} as an exact sequence consisting of projective
objects. The Proposition is proved. 
\end{proof}


This proposition implies that every projective object in a finite tensor 
category is injective, and vice versa. This is a generalization of the fact
that a finite dimensional Hopf algebra is Frobenius. It also implies that 
an indecomposable projective object $P$ has a unique simple subobject.
This subobject will be called the {\it socle} of $P$. 

\subsection{Surjective quasi-tensor functors}
Let $\mC$, $\mD$ be abelian categories. 
Let $F: \mC\to \mD$ be an additive functor. 

\begin{definition} We will say that $F$ is surjective if 
any object of $\mD$ is a subquotient in 
$F(X)$ for some $X\in \mC$. 
\end{definition}

{\bf Example.} Let $A,B$ be coalgebras, and $f:A\to B$ 
a homomorphism. Let $F=f^*: A-\mbox{comod}\to B-\mbox{comod}$ be the 
corestriction functor.
Then $F$ is surjective iff $f$ is surjective. 

Now let $\mC$, $\mD$ be finite tensor categories.
An additive functor $F: \mC\to \mD$ is said to be quasi-tensor if 
it is exact and faithful, and for any 
objects $X,Y$, the object $F(X)\otimes F(Y)$ is isomorphic to 
$F(X\otimes Y)$. In particular a tensor functor is quasi-tensor. 

\begin{theorem} \label{projsurj}
Let $F: \mC\to \mD$ be a surjective quasi-tensor functor. Then 
$F$ maps projective objects to projective ones. 
\end{theorem}

The proof of this theorem is given later in this section. 

\subsection{Frobenius-Perron dimensions}\label{FPdim}

Let $\mC$ be a finite tensor category. 
Recall \cite{E} that for each object $X$ of $\mC$ one can define its 
Frobenius-Perron 
dimension $d_+(X)$, which is additive on exact sequences and multiplicative
(namely, $d_+(X)$ is the largest positive eigenvalue of the matrix of 
left or right multiplication by $X$). This is an algebraic integer. 
The function $d_+$ is the unique character of $Gr(\mC)$ which takes positive 
values on all simple objects $X$ of $\mC$
(this follows from the Frobenius-Perron theorem). 
Therefore, any quasi-tensor functor between finite tensor categories preserves 
Frobenius-Perron dimensions. 

The following statement is well known in the semisimple case, see e.g. 
\cite{ENO}.

\begin{proposition}\label{qua} 
Assume that the Frobenius-Perron dimensions of objects
in $\cC$ are integers. Then $\cC$ is equivalent to the representation
category of a finite dimensional quasi-Hopf algebra. \end{proposition}

\begin{proof} Clearly it is enough to construct an exact functor 
$F: \cC \to \Ve_k$ together with a functorial isomorphism $F(X\otimes Y)\simeq
F(X)\otimes F(Y)$. Define $P=\oplus_iP_i^{d_+(L_i)}$ and $F=\Hom(P,?)$. 
Obviously, $F$ is exact and $\dim F(X)=d_+(X)$. Using \cite{De2} Proposition
5.13 (vi) we continue the functors $F(?\otimes ?)$ and $F(?)\otimes F(?)$ to
the functors $\cC \boxtimes \cC \to \Ve_k$. Both of these functors are exact 
and take the same values on the simple objects of $\cC \boxtimes \cC$. Thus
these functors are isomorphic and we are done. \end{proof}

Let $A$ be a separable algebra. For a finite tensor category $\cC$ and a
tensor functor $F: \cC \to \Bimod(A)$ one constructs a tensor equivalence of
$\cC$ and the representation category of a finite dimensional weak Hopf algebra
(see e.g. \cite{ENO}). If the functor $F$ is assumed to be only quasi-tensor
one should replace weak Hopf algebras by weak quasi-Hopf algebras, see
\cite{MaS}. We have the following

\begin{proposition} Any finite tensor category $\cC$ is equivalent to the 
representation category of a finite dimensional weak quasi-Hopf algebra. 
\end{proposition}

\begin{proof} We need to construct a quasi-tensor functor $F: \cC \to 
\Bimod(A)$. Set $A=\oplus_{i\in I}ke_i$, $e_ie_j=\delta_{ij}e_i$. Let $A_{ij}$
denote the $A-$bimodule $e_iAe_j$. Set $F(X)=\oplus_{i,j\in I}
\Hom(P_i, X\otimes L_j)\otimes A_{ij}$. By the same argument as in the
proof of Proposition \ref{qua} the functors $F(X\otimes Y)$ and $F(X)\otimes_A
F(Y)$ are isomorphic. The Proposition is proved. \end{proof}

\subsection{Projectivity defect}

Let $\mC$ be a finite tensor category, and $X\in \mC$.
Let us write $X$ as a direct sum of indecomposable objects (such a representation 
is unique).  
Define the projectivity defect $p(X)$ of $X$ to be the sum of 
Frobenius-Perron dimensions of all the non-projective summands in this sum. 
It is clear that $p(X\oplus Y)= p(X)+p(Y)$. Also, it follows from 
Proposition \ref{tensprod} that $p(X\otimes Y)\le p(X)p(Y)$.  

\subsection{Proof of Theorem \ref{projsurj}}

Let $P_i$ be the indecomposable projective objects in 
$\mC$. Let $P_i\otimes P_j\cong \oplus_k B_{ij}^kP_k$, and let 
$B_i$ be the matrix with entries $B_{ij}^k$. Also, let $B=\sum
B_i$. Obviously, $B$ has strictly positive entries, and 
the Frobenius-Perron eigenvalue of $B$ is $\sum_i d_+(P_i)$.

On the other hand, let $F:\mC\to \mD$ be a surjective quasi-tensor functor
between finite tensor categories. Let $p_j=p(F(P_j))$, and $\bold
p$ be the vector with entries $p_j$. 
Then we get $p_ip_j\ge \sum_k B_{ij}^kp_k$, 
so $(\sum_i p_i)\bold p\ge B\bold p$. So, 
either $p_i$ are all 
zero, or they are all positive, and the norm of $B$ 
with respect to the norm $|x|=\sum p_i|x_i|$ is at most $\sum
p_i$. Then, by the Frobenius-Perron theorem, one would have 
$p_i=d_+(P_i)$ for all $i$.

 Assume the second option is the case. 
Then $F(P_i)$ do not contain nonzero 
projective objects as direct summands, and hence 
for any projective $P\in \mC$, $F(P)$ cannot contain a nonzero projective 
object as a direct summand. However, let $Q$ be a projective object  
 of $\mD$. Then there exists an object $X\in \mC$ such that 
$Q$ is a subquotient of $F(X)$. Since any $X$ is a quotient of a projective 
object, and $F$ is exact, we may assume that $X=P$ is projective.  
So $Q$ occurs as a subquotient in $F(P)$. As $Q$ is both projective and injective, 
it is actually a direct summand in $F(P)$. Contradiction. 

Thus, $p_i=0$ and $F(P_i)$ are projective. 
The theorem is proved. 

\subsection{Categorical freeness}\label{freen} 

Let $K(\mC)$ denote the free abelian group generated by 
isomorphism classes of indecomposable projective objects of
$\mC$. Elements of 
$K(\mC)\otimes \Bbb C$ will be called virtual projective
objects. Recall \cite{E} that for every finite tensor category $\mC$ one may 
define a unique virtual projective 
object $R_\mC:=\sum d_+(L_i)P_i\in 
K(\mC)\otimes_{\Bbb Z}\Bbb C$, 
such that $X\otimes R_\mC=R_\mC\otimes X=d_+(X)R_\mC$ for any
$X\in \mC$,
and $\dim \Hom(R_\mC,\bold 1)=1$. 
The Frobenius-Perron dimension of this object is called the Frobenius-Perron 
dimension of $\mC$, and denoted $d_+(\mC)$.  

{\bf Remark.} We note the following useful inequality: 
$d_+(\mC)\ge Nd_+(P)$, where $N$ is the number of simple 
objects in $\mC$, and $P$ is the projective cover of the neutral
object $\bold 1$. Indeed, for any simple object $V$ the projective object 
$P(V)\otimes {}^*V$ has a nontrivial homomorphism to $\bold 1$,
and hence contains $P$. So $d_+(P(V))d_+(V)\ge d_+(P)$. 
Adding these inequalities over all simple $V$, we get the
result. 

For any surjective quasi-tensor functor $F:\mC\to \mD$, one has 
$$F(R_\mC)=\frac{d_+(\mC)}{d_+(\mD)} R_\mD. \eqno(1)$$ 
Indeed, by Theorem \ref{projsurj}, $F(R_\mC)$ is a virtually
projective object. Thus, $F(R_\mC)$ must be proportional to $R_\mD$, since both 
(when written in the basis $P_i$) are 
eigenvectors of a matrix with strictly positive entries
with its Frobenius-Perron eigenvalue. (For this matrix we may take 
the matrix of multiplication by $F(X)$, 
where $X$ is such that $F(X)$ contains as constituents all 
simple objects of $\mD$; such exists by the surjectivity of $F$).
This shows that $d_+(\mC)\ge d_+(\mD)$, and $d_+(\mD)$ divides $d_+(\mC)$ 
as an algebraic integer: in fact,  
$\frac{d_+(\mC)}{d_+(\mD)}=\sum d_+(L_i)\dim\Hom(F(P_i),\bold 1_\mD)$. 

Suppose now that the Frobenius-Perron dimensions of objects in $\mC$ are integers 
and thus $\mC$ is the representation category of a quasi-Hopf algebra, see
Proposition \ref{qua}. In this case $R_\mC$ is an honest
(not only virtual) projective object of $\mC$. Multiples of $R_\mC$ 
will be called free objects of $\mC$, and the multiplicity 
will be refereed to as rank.  

Then Theorem \ref{projsurj} and the fact that $F(R_{\mC})$ is
proportional to $R_\mD$ implies 

\begin{corollary} (The categorical freeness theorem) 
If the Frobenius-Perron dimensions in $\mC$ are integers, 
and $F:\mC\to \mD$ is a surjective 
quasi-tensor functor then the Frobenius-Perron dimensions in
$\mD$ are integers as well, and the object $F(R_{\mC})$  
is free of rank $d_+(\mC)/d_+(\mD)$ (which is an integer).  
\end{corollary}

In the Hopf case this theorem is well known
and much used; it is due to Nichols and Zoeller \cite{NZ}, 
and claims that a finite dimensional Hopf algebra is free as 
a module over a Hopf subalgebra. 
In the quasi-Hopf case 
it was recently proved in \cite{ENO} 
in the semisimple case, and in general 
by Schauenburg \cite{Sch}, Theorem 3.2.

\subsection{The distinguished character}

Since duals to projective objects are projective, 
we can define a map $D: I\to I$ such that $P_i^*=P_{D(i)}$. 
It is clear that $D^2(i)=i^{**}$. 

Let $0$ be the label for the identity object. Let $\rho=D(0)$.
We have $\Hom(P_i^*,L_j)=\Hom(\bold 1,P_i\otimes L_j)=
\Hom(\bold 1,\oplus_k N_{kj^*}^iP_k)$. This space has dimension 
$N_{\rho j^*}^i$. Thus we get 
$$
N_{\rho j^*}^i=\delta_{D(i),j}.
$$
Let now $L_\rho$ be the corresponding simple object. 
We have 
$$
L_\rho^*\otimes P_m\cong \oplus_k N_{\rho m}^k P_k\cong P_{D(m)^*}.
$$
 
\begin{lemma}
$L_\rho$ is an invertible object.
\end{lemma}

\begin{proof}
The last equation implies that the matrix of action of
$L_{\rho^*}$ on projectives is a permutation matrix. 
Hence, the Frobenius-Perron dimension of $L_{\rho^*}$ is $1$, and
we are done. 
\end{proof}

\begin{lemma} \label{duali}
One has: $P_{D(i)}=P_{{}^*i}\otimes L_\rho$;
$L_{D(i)}=L_{{}^*i}\otimes L_\rho$.
\end{lemma}

\begin{proof}
It suffices to prove the first statement. 
Therefore we need to show that 
$\dim \Hom(P_i^*,L_j)=\dim\Hom(P_{{}^*i},L_j\otimes L_{\rho^*})$. 
The left hand side was computed before, it is 
$N_{\rho j^*}^{i}$. On the other hand, the right hand side is 
$N_{j,\rho^*}^{{}^*i}$ (we use that $\rho^*=^*\rho$ for an invertible 
object $\rho$). 
These numbers are equal by the definition of
duality, so we are done. 
\end{proof}

\begin{corollary} One has: 
$P_{i^{**}}=L_\rho^*\otimes P_{{}^{**}i}\otimes L_\rho$;
$L_{i^{**}}=L_\rho^*\otimes L_{{}^{**}i}\otimes L_\rho$.
\end{corollary}

\begin{proof} Again, it suffices to prove the first statement. 
We have 
$$
P_{i^{**}}=P_i^{**}=(P_{{}^*i}\otimes L_\rho)^*=
L_\rho^*\otimes P_{{}^*i}^*=L_\rho^*\otimes P_{{}^{**}i}\otimes L_\rho
$$
\end{proof} 

\begin{definition} $L_\rho$ is called the distinguished invertible object 
of $\mC$.
\end{definition}

\begin{proposition}
Let $H$ be a finite dimensional Hopf algebra, and $\mC=\Rep(H)$. 
Then $L_\rho$ is the distinguished group-like element of $H^*$.
\end{proposition}

\begin{proof}
Let $\chi$ be the distinguished character of $H$. 
Then there exists a nonzero element 
$I\in H$ such that $xI=\varepsilon(x)I$
(i.e. $I$ is a left integral) and 
$Ix=\chi(x)I$. This means that for any 
$V\in \mC$, $I$ defines a morphism from $V\otimes \chi^{-1}$ 
to $V$. 

The element $I$ belongs to the submodule $P_i$ of 
$H$, whose socle (i.e. the irreducible submodule) 
is the trivial $H$-module. Thus, $P_i^*=P_1$, and hence
by Lemma \ref{duali}, $i=\rho$. Thus, 
$I$ defines a nonzero 
(but rank $1$) 
morphism $P_\rho\otimes \chi^{-1}\to P_\rho$.
The image of this morphism, because of rank $1$, must be $L_0=\bold 1$, 
so $\bold 1$ is a quotient of $P_\rho\otimes \chi^{-1}$, and hence 
$\chi$ is a quotient of $P_\rho$. Thus, $\chi=L_\rho$, and we are done. 
\end{proof}

\begin{remark}
A similar proof applies to weak Hopf algebras.
\end{remark}

\begin{conjecture} For any finite tensor category $\mC$, 
there exists a natural isomorphism of tensor functors
$V^{**}\to L_\rho^*\otimes ^{**}V\otimes L_\rho$.
\end{conjecture}

For Hopf algebras, this follows from Radford's formula for $S^4$. 
For weak Hopf algebras, it follows from the Nikshych's generalization of 
Radford's formula, see \cite{N}, \cite{ENO}.

\subsection{Dimensions of projective objects and degeneracy 
of the Cartan matrix}

The following result in the Hopf algebra case was proved by M.Lorenz \cite{L}; 
our proof in the categorical setting is analogous to his. 

Let $C_{ij}=[P_i:L_j]$ be the entries of the Cartan matrix of $\mC$. 

\begin{theorem}\label{lorentz}
Suppose that $\mC$ is not semisimple, and admits 
an isomorphism 
of additive functors $u: {\rm Id}\to **$
(for example, $\mC$ is braided).
Then the Cartan matrix $C$ is degenerate over the ground field $k$.
\end{theorem}

\begin{proof}
Let $\dim(V)=\Tr|_V(u)$ be the dimension function 
defined by the categorical trace of $u$. 
Then the dimension of every projective object $P$ is zero.
Indeed, the dimension of $P$ is 
the composition of maps 
$\bold 1\to P\otimes P^*\to P^{**}\otimes P^*\to \bold 1$, where 
the maps are the coevaluation, $u\otimes 1$, and the evaluation. 
If this map is nonzero then $\bold 1$ is a direct summand in 
$P\otimes P^*$, which is projective. Thus $\bold 1$ is projective,
hence any object $V=V\otimes \bold 1$ is projective. 
So $\mC$ is semisimple. Contradiction.

Since the dimension of the trivial object $\bold 1$ cannot be zero, 
$\bold 1$ is not a linear combination of projective objects 
in the Grothendieck group tensored with $k$. We are done. 
\end{proof}

\subsection{Absence of primitive elements}

The following theorem is a categorical version 
of the absence of primitive elements in finite dimensional Hopf algebras
in characteristic zero. Again, the proof is a categorical version of the 
standard proof for Hopf algebras. 

\begin{theorem} \label{abs} Assume that $k$ has characteristic 0.
Let $\mC$ be a finite tensor category over $k$.
Then ${\rm Ext}^1(\bold 1,\bold 1)=0$. 
\end{theorem}

\begin{proof}
Assume the contrary, and suppose 
that $V$ is a nontrivial extension of $\bold 1$ by itself. 
Let $P$ be the projective cover of $\bold 1$.
Then $\Hom(P,V)$ is a 2-dimensional space, with a filtration 
induced by the filtration on $V$. 
Let $v_0,v_1$ be a basis compatible to the filtration, 
i.e. $v_0$ spans the 1-dimensional subspace 
defined by the filtration.
Let $A=\End(P)$ (this is a finite dimensional algebra). 
Let $\varepsilon: A\to \Bbb C$ be the character 
defined by the action of $A$ on $\Hom(P,\bold 1)$. 
Then the matrix of $a\in A$ in $v_0,v_1$ has the form
\begin{equation}
[a]_1=\begin{pmatrix}
\varepsilon(a) & \chi_1(a) \\
0 & \varepsilon(a)
\end{pmatrix}
\end{equation}
where $\chi_1\in A^*$ is nonzero. 
Since $a\to [a]_1$ is a homomorphism, $\chi_1$ is a derivation:
$\chi_1(xy)=\chi_1(x)\varepsilon(y)+\varepsilon(x)\chi_1(y)$.

Now consider the representation $V\otimes V$. 
The space $\Hom(P,V\otimes V)$ 
is 4-dimensional, and has a 3-step filtration, 
with basis $v_{00}; v_{01},v_{10}; v_{11}$, consistent with this filtration. 
The matrix of $a\in \End(P)$ in this basis (under appropriate normalization 
of basis vectors) is 

\begin{equation}
[a]_2=\begin{pmatrix}
\varepsilon(a) & \chi_1(a) & \chi_1(a) & \chi_2(a)  \\
0 & \varepsilon(a) & 0 & \chi_1(a)  \\
0 & 0 & \varepsilon(a) & \chi_1(a)  \\
0 & 0 & 0 & \varepsilon(a)
\end{pmatrix}
\end{equation}
Since $a\to [a]_2$ is a homomorphism, 
we find 
$$
\chi_2(ab)=\varepsilon(a)\chi_2(b)+\chi_2(a)\varepsilon(b)+2\chi_1(a)\chi_1(b).
$$
We can now proceed further (i.e. consider $V\otimes V\otimes V$ etc.)
and define for every positive $n$, 
a linear function $\chi_n\in A^*$ 
which satisfies the equation 
$$
\chi_n(ab)=\sum_{j=0}^n\begin{pmatrix} n \\ j \end{pmatrix}
\chi_j(a)\chi_{n-j}(b),
$$
where $\chi_0=\varepsilon$. 
Thus for any $s\in k$, we can define 
$\phi_s: A\to k((t))$ by $\phi_s(a)=\sum_{m\ge 0}\chi_m(a)s^mt^m/m!$, 
and we find that $\phi_s$ is a family of pairwise distinct homomorphisms. 
This is a contradiction, as $A$ is a finite dimensional algebra. 
We are done.
\end{proof}
In particular, if $\cC$ has a unique simple object then $\cC$ is equivalent
to the category $\Ve_k$ of vector spaces. Certainly this is not true in
characteristic $p>0$, a counterexample being $\cC =\Rep(G)$ for a finite
$p-$group $G$.

\subsection{The Ext algebra of a finite tensor category} 
We expect the following to be true.

\begin{conjecture} For a finite tensor category $\cC$ the algebra 
${\rm Ext}^*(\be,\be)$ is finitely generated. Moreover, for any $X\in \cC$ the
module ${\rm Ext}^*(\be,X)$ over ${\rm Ext}^*(\be,\be)$ is finitely generated.
\end{conjecture}

Note that the algebra ${\rm Ext}^*(\be,\be)$ is graded commutative, see
e.g. \cite{SA} and references therein. 

It is known that the conjecture is true for $\cC =\Rep(H)$ where $H$ is
either commutative or cocommutative Hopf algebra (the first since the algebra
of functions on a finite group scheme is a complete intersection and the
second is a deep Theorem of E.~Friedlander and A.~Suslin \cite{FS}).

\subsection{Surjective and injective functors}
A tensor subcategory in a finite tensor category $\cD$ 
is a full subcategory $\cC\subset \cD$ 
which is closed under taking subquotients, tensor products, 
duality, and contains
the neutral object. 
A tensor functor $F: \cC\to \cD$ is injective if it is an equivalence of 
$\cC$ onto a tensor subcategory of $\cD$. 

\begin{proposition}\label{subc}
If $F: \cC\to \cD$ is injective, then 
$d_+(\cC)\le d_+(\cD)$. The equality is achieved 
if and only if $F$ is an equivalence. 
\end{proposition}

\begin{proof} We may assume that $\cC$ is a tensor subcategory 
of $\cD$. Let $X$ be a simple object of $\cC$. 
Then $X$ is also a simple object of $\cD$. 
Let $P_{\cC}(X),P_{\cD}(X)$ be the projective covers
of $X$ in $\cC$, $\cD$. Then we have 
projections $a_{\cC}: P_{\cC}(X)\to X$, 
$a_{\cD}: P_{\cD}(X)\to X$ (in $\cC,\cD$, respectively),
and there exists a morphism 
$b: P_{\cD}(X)\to P_{\cC}(X)$ such that 
$a_{\cC}\circ b=a_{\cD}$. We claim that $b$ is onto. 
Indeed, assume the contrary. Let $L$ be a simple quotient 
of the cokernel of $b$, and $f$ 
a projection $P_{\mC}(X)\to L$. 
 It is clear that $L=X$ and $f$ is proportional to $a_{\mC}$.
But this is a contradiction, since $f\circ b=0$. 

Thus, $d_+(P_{\cD}(X))\ge d_+(P_{\cC}(X))$. 
This implies the first statement of the proposition. 

Let us now prove the second statement. 
The equality $d_+(\cC)=d_+(\cD)$ implies
that 1) all simple objects of $\cC$ are also contained in $\cD$, 
and 2) $b$ is an isomorphism, i.e. $P_{\cC}(X)=P_{\cD}(X)$
(i.e., $\cC$ is a Serre subcategory of $\cD$). This implies that 
$\cC=\cD$. 
\end{proof}

\begin{proposition}\label{samedim} 
If $F: \cC\to \cD$ is surjective, then 
$d_+(\cC)\ge d_+(\cD)$. The equality is achieved 
if and only if $F$ is an equivalence. 
\end{proposition}

\begin{proof} 
The first statement has already been proved, so let us 
prove the second one. 
Let $L_i,P_i, 1\le i\le n$ be the simple and projective objects 
in $\cC$ and $L_j',P_j', 1\le j\le m$ the simple and projective objects in
$\cD$. Let $d_i,d_j'$ be the dimensions of $L_i,L_j'$.
We have $F(L_i)=\sum_j a_{ij}L_j'$ (in the Grothendieck group
$Gr(\mC)$), and
$F(P_i)=\oplus_j b_{ji}P_j'$ (in $K(\mC)$). From the first equation we get
$d_i=\sum a_{ij}d_j'$, and since $F(R_{\cC})=
R_{\cD}$, from the second equation we get $\sum b_{ji}d_i=d_j'$, 
So if $\bold d$, $\bold d'$ are the  vectors with 
entries $d_i,d_j'$, and $A,B$ matrices with entries $a_{ij}$,
$b_{ji}$, then $A\bold d'=\bold d$, $B\bold d=\bold d'$. Thus 
$AB\bold d=\bold d$. 

Now, the functor $F$ is faithful. Hence, ${\rm
Hom}(F(P_i),F(L_i))\ne 0$. But the dimension of this space is 
$\sum_j a_{ij}b_{ji}=(AB)_{ii}$. Hence the diagonal entries 
of $AB$ are $\ge 1$. Since $AB\bold d=\bold d$,
the entries of $\bold d$ are positive, and 
the entries of $AB$ are nonnegative, we conclude that
$AB=1$. 

We will now show that $n\ge m$. This will imply that $BA=1$. 
Since $AB=1$, 
for any $i$ there exists a unique $j$ such that 
$a_{ij}b_{ji}\ne 0$; call it $j(i)$.
It suffices to show that for any $j$ there exists $i$ such that 
$j=j(i)$. Assume the contrary, i.e. some $j\ne j(i)$ for any $i$. 
Then $a_{ij}b_{jk}=0$ for all $i,k$. 
Choose $i$ so that $L_j'$ is contained as a constituent in
$F(L_i)$ (it must exist as $F$ is surjective).   
Then $a_{ij}\ne 0$, so $b_{jk}=0$ for all $k$. 
This means that $P_j'$ is not a direct summand 
of $F(P_k)$ for any $k$, i.e. 
is not a subquotient of $F(Q)$ 
for any projective object $Q$. Contradiction with 
surjectivity of
$F$. 

Thus $AB=1$, $BA=1$. This means that $A$ is a permutation
matrix, and $B=A^{-1}$. This easily implies that $F$ is an
equivalence. 
\end{proof} 

Let $F:\cC\to \cD$ be a tensor functor between two finite tensor categories. 
Then we can define a tensor subcategory ${\rm Im}F$ 
of $\cD$ to be the full subcategory of $\cD$ 
consisting of objects contained as subquotients in $F(X)$ for
some $X\in \cC$. 
It is clear that ${\rm Im}F$ is a tensor subcategory
of $\cD$. 

The functor $F$ is naturally written as a composition of two 
tensor functors: $F=F_i\circ F_s$, where $F_s: \cC\to {\rm Im}F$ 
is surjective, and $F_i: {\rm Im}F\to \cD$ 
is injective. Clearly, $F$ is surjective iff 
$F_i$ is an equivalence, and $F$ is injective iff $F_s$ 
is an equivalence. 

\begin{corollary} \label{injsurj}
(i) If $d_+({\rm Im}F)=d_+(\cC)$ iff $F$ is injective.

(ii) If $d_+({\rm Im}F)=d_+(\cD)$ iff $F$ is surjective.
\end{corollary}

\begin{proof}
(i) Follows from Proposition \ref{subc}.

(ii) Follows from Proposition \ref{samedim}. 
\end{proof}

\begin{corollary} \label{notinjsurj}
Suppose that a tensor functor $F:\cC\to \cD$ 
between finite tensor categories
factors through a finite tensor category 
${\mathcal E}$, such that $d_+({\mathcal E})<{\rm min}(d_+(\cC),d_+(\cD))$. 
Then $F$ is neither surjective nor 
injective.  
\end{corollary}

\begin{proof} We have $F=F_1\circ F_2$, $F_2: \cC\to {\mathcal
E}$, $F_1: {\mathcal E}\to \cD$. Clearly, ${\rm Im}F$ is a tensor
subcategory in ${\rm
Im}F_1$, so $d_+({\rm Im} F)\le 
d_+({\rm Im}F_1)\le d_+({\mathcal E})$. Thus, by Proposition 
\ref{injsurj}, $F$ is neither surjective nor injective. 
\end{proof} 

\section{Exact module categories} In this section we will work with slightly
more general categories than finite tensor categories. Namely by a multi-tensor
category we will mean a rigid tensor category in which the unit object $\be$
is completely reducible, $\be =\oplus_{i\in I}\be_i$. A multi-tensor category
$\cC$ is called indecomposable if the subcategory
$\be_i\otimes \cC \otimes \be_j\subset \cC$ is nonzero for all $i,j\in I$.
Note that $\cC \simeq \oplus_{i,j\in I}\be_i\otimes \cC\otimes \be_j$. 
In what follows all multi-tensor categories are assumed to be indecomposable.
We leave to the reader to check that Propositions \ref{tensprod}, 
\ref{dualproj} and Theorems \ref{projsurj}, \ref{abs} remain true for
finite multi-tensor categories.
 
\subsection{Definition and basic properties} We will assume in what follows
that categories $\cM, \cM_1, \cM_2$ etc have only finitely many isomorphism
classes of simple objects.

\begin{definition} Let $\cC$ be a finite multi-tensor category. A module 
category
$\cM$ over $\cC$ is called exact if for any {\em projective} object $P\in \cC$
and any object $X\in \cM$ the object $P\otimes X\in \cM$ is projective.
\end{definition}

Our aim is to show that the notion of an exact module category is a good 
generalization of the notion of a semisimple module category over a fusion
category.

\begin{remark} (i) Let $\cM$ be an arbitrary module category over $\cC$. For 
any projective object $Q\in \cM$ and any object $L\in \cC$ the object 
$L\otimes Q$ is projective. Indeed, the functor $\Hom(L\otimes Q,?)$ is 
isomorphic to $\Hom(Q,{}^*L\otimes ?)$ and hence is exact.

(ii) We will show later (Proposition \ref{flat}) that any module functor from
an exact module category is exact. This explains our choice for this name. 
\end{remark}

\begin{example}\label{exex}
 (i) Any finite tensor category $\cC$ considered as a module 
category over itself is exact. Also, the category $\cC$ considered as a module
category over $\cC \boxtimes \cC^{op}$ is exact (here $\cC^{op}$ is the same
category as $\cC$ but with new tensor product $X\tilde \otimes Y:=Y\otimes X$;
$\cC$ is a module category over $\cC \boxtimes \cC^{op}$ via $(X\boxtimes Y)
\otimes Z=X\otimes Z\otimes Y$).

(ii) Let $F: \cC \to \cD$ be a surjective tensor functor. Then the category
$\cD$ considered as a module category over $\cC$ is exact by Theorem 
\ref{projsurj}.

(iii) Assume that $\cC$ is a semisimple category (thus $\cC$ is a fusion 
category). A module category over $\cC$ is exact if and only if it is
semisimple. Indeed, in this case the unit object of $\cC$ is projective.
\end{example}

\begin{lemma} Let $\cM$ be an exact module category over $\cC$. The category
$\cM$ has enough projective objects. In particular the category $\cM$ is
finite.
\end{lemma}

\begin{proof} Let $P_0\in \cC$ denote the projective cover of the unit object
in $\cC$. Then the natural map $P_0\otimes X\to X$ is surjective for any $X\in
\cM$ and $P_0\otimes X$ is projective by definition of an exact module 
category. \end{proof}

\begin{lemma}\label{proinj} Let $\cM$ be an exact module category over $\cC$.
Let $P\in \cC$ be projective and $X\in \cM$. Then $P\otimes X$
is injective.
\end{lemma}

\begin{proof} The functor $\Hom(?,P\otimes X)$ is isomorphic to the functor
$\Hom(P^*\otimes ?,X)$. The object $P^*$ is projective by Proposition 
\ref{dualproj}. Thus
for any exact sequence \linebreak
$0\to Y_1\to Y_2\to Y_3\to 0$ the sequence
$0\to P^*\otimes Y_1\to P^*\otimes Y_2\to P^*\otimes Y_3\to 0$ splits and
hence the functor $\Hom(P^*\otimes ?,X)$ is exact. The Lemma is proved.
\end{proof}

\begin{corollary} In the category $\cM$ any projective object is injective and
vice versa.\end{corollary}

\begin{proof} Any projective object of $\cM$ is a direct summand of the object
of the form $P_0\otimes X$ and thus is injective. \end{proof}

\begin{remark} A finite abelian category $\cA$ is called a Frobenius category
if any projective object of $\cA$ is injective and vice versa. Thus any
exact module category over a finite multi-tensor category (in particular any
finite multi-tensor category itself) is a Frobenius category. It is well known 
that any object of a Frobenius category admitting a finite projective 
resolution is projective (indeed, the last nonzero arrow of this resolution is
an embedding of projective ($=$ injective) modules and therefore is an 
inclusion of a direct summand. Hence the resolution can be replaced by a 
shorter one and by induction we are done). Thus any Frobenius category is 
either semisimple or of infinite homological dimension. \end{remark}

Let $Irr(\cM)$ denote the set of (isomorphism classes of) simple objects in
$\cM$. Let us introduce the following relation on $Irr(\cM)$:
two objects $X,Y\in Irr(\cM)$ are related if $Y$ appears as a subquotient
of $L\otimes X$ for some $L\in \cC$. 

\begin{lemma} The relation above is reflexive, symmetric and
transitive.\end{lemma}

\begin{proof} Since $\be \otimes X=X$ we have the reflexivity. Let $X, Y, Z\in 
Irr(\cM)$ and $L_1, L_2\in \cC$. If $Y$ is a subquotient of $L_1\otimes X$ and
$Z$ is a subquotient of $L_2\otimes Y$ then $Z$ is a subquotient of 
$(L_2\otimes L_1)\otimes X$ (since $\otimes$ is exact) whence we get the 
transitivity. Now assume that $Y$ is a subquotient of $L\otimes X$. Then the 
projective cover $P(Y)$ of $Y$ is a direct summand of $P_0\otimes L\otimes X$;
hence there exists $S\in \cC$ such that $\Hom(S\otimes X,Y)\ne 0$ (for example
$S=P_0\otimes L$). Thus $\Hom(X,S^*\otimes Y)=\Hom(S\otimes X,Y)\ne 0$ and
hence $X$ is a submodule of $S^*\otimes Y$. Consequently our equivalence
relation is symmetric. \end{proof}

Thus our relation is an equivalence relation. Hence $Irr(\cM)$ is 
partitioned into equivalence classes, $Irr(\cM)=\bigsqcup_{i\in I}
Irr(\cM)_i$. For an equivalence class $i\in I$ let $\cM_i$ denote the full 
subcategory of $\cM$ consisting of objects all simple subquotients of which 
lie in $Irr(\cM)_i$. Clearly, $\cM_i$ is a module subcategory of $\cM$.

\begin{proposition} The module categories $\cM_i$ are exact. The category 
$\cM$  is the direct sum of its module subcategories $\cM_i$. 
\end{proposition}

\begin{proof} For any $X\in Irr(\cM)_i$ its projective cover is a direct 
summand of $P_0\otimes X$ and hence lies in the category $\cM_i$. Hence the 
category $\cM$ is the direct sum of its subcategories $\cM_i$, and $\cM_i$
are exact.
\end{proof}

Recall (see e.g. \cite{O}) that a $\BZ_+-$module over a $\BZ_+-$ring is 
called irreducible if it has no notrivial $\BZ_+-$submodules. 

\begin{corollary} Let $\cM$ be an indecomposable exact module category. 
Then the Grothendieck group $Gr(\cM)$ is an irreducible $\BZ_+-$module over
$Gr(\cC)$. \end{corollary}

In particular, for a given category $\cC$ there are only finitely many
$\BZ_+-$modules over $Gr(\cC)$ which are of the form $Gr(\cM)$ where $\cM$ is
an indecomposable exact module category over $\cC$, see \cite{O}, Proposition 
2.1.

The crucial property of exact module categories is the following

\begin{proposition} \label{flat}
Let $\cM_1$ and $\cM_2$ be two module categories over $\cC$. Assume that
$\cM_1$ is exact. Then any additive module functor $F:\cM_1\to \cM_2$ is exact.
\end{proposition}

\begin{proof} Let $0\to X\to Y\to Z\to 0$ be an exact sequence in $\cM_1$.
Assume that the sequence $0\to F(X)\to F(Y)\to F(Z)\to 0$ is not exact. Then
the sequence $0\to P\otimes F(X)\to P\otimes F(Y)\to P\otimes F(Z)\to 0$ is
also nonexact for any nonzero object $P\in \cC$ since the functor $P\otimes ?$
is exact and $P\otimes X=0$ implies $X=0$. In particular we can take $P$ to be
projective. But then the sequence $0\to P\otimes X\to P\otimes Y\to 
P\otimes Z\to 0$ is exact and split and hence the sequence $0\to 
F(P\otimes X)\to F(P\otimes Y)\to F(P\otimes Z)\to 0$ is exact and we 
get a contradiction. \end{proof}

\begin{remark}
We will see later that this Proposition actually characterizes exact
module categories.\end{remark}

\subsection{Morita theory} An important technical tool in the study of
module categories is the notion of internal Hom. Let $\cM$ be a module
category over $\cC$ and $M_1, M_2\in \cM$. Consider the functor 
$\Hom(?\otimes M_1,M_2)$ from the category $\cC$ to the category of vector
spaces. This functor is left exact and thus is representable (see e.g.
\cite{Se}, Chap. II, \S 4). 

\begin{definition} The internal Hom $\iHom(M_1,M_2)$ is an object of $\cC$
representing the functor $\Hom(?\otimes M_1,M_2)$.\end{definition}

Note that by Yoneda's Lemma $\iHom(M_1,M_2)$ is a bifunctor.

\begin{lemma} \label{can} There are canonical isomorphims
$$\begin{array}{cl}
(1)&\Hom(X\otimes M_1,M_2)\cong \Hom(X,\iHom(M_1,M_2)),\\
(2)&\Hom(M_1,X\otimes M_2)\cong \Hom(\be,X\otimes \iHom(M_1,M_2)),\\
(3)&\iHom(X\otimes M_1,M_2)\cong \iHom(M_1,M_2)\otimes X^*,\\
(4)&\iHom(M_1,X\otimes M_2)\cong X\otimes \iHom(M_1,M_2).\end{array}$$
\end{lemma}

\begin{proof} See \cite{O}, Lemma 3.3. \end{proof}

Note that isomorphisms (3) and (4) above show that $\iHom(?,?)$ is a module
functor in each variable. Thus by Proposition \ref{flat} we have

\begin{corollary}\label{biex}
 Assume that $\cM$ is exact module category. Then the
functor $\iHom(?,?)$ is biexact. \end{corollary} 

The mere definition of the internal Hom allow us to prove the converse to
Proposition \ref{flat}:

\begin{proposition} Let $\cM_1, \cM_2$ be two nonzero module categories over
$\cC$. Assume that any module functor from $\cM_1$ to $\cM_2$ is exact. Then
the module category $\cM_1$ is exact.\end{proposition}

\begin{proof} First we claim that under our assumptions any module functor
$F\in Fun_{\cC}(\cM_1,\cC)$ is exact. Indeed let $0\ne M\in \cM_2$. The
functor $F(?)\otimes M\in Fun_{\cC}(\cM_1,\cM_2)$ is exact. Since $?\otimes M$
is exact and $X\otimes M=0$ implies $X=0$ we see that $F$ is exact.

In particular for any object $N\in \cM_1$ the functor $\iHom(N,?):\cM_1\to 
\cC$ is exact since it is a module functor. Now let $P\in \cC$ be any 
projective object. Then for any $N\in \cM_1$ one has $\Hom(P\otimes N,?)=
\Hom(P,\iHom(N,?))$ and thus the functor $\Hom(P\otimes N,?)$ is exact. By 
the definition of an exact module category we are done. \end{proof}

For two objects $M_1, M_2$ of a module category $\cM$ we have the canonical 
morphism 
$$ev_{M_1,M_2}: \iHom(M_1,M_2)\otimes M_1\to M_2$$ 
obtained as the image of $id$ under the isomorphism
$$\Hom(\iHom(M_1,M_2),\iHom(M_1,M_2))\cong 
\Hom(\iHom(M_1,M_2)\otimes M_1,M_2).$$
Let $M_1, M_2, M_3$ be three objects of $\cM$. Then there is a canonical
composition morphism
$$(\iHom(M_2,M_3)\otimes \iHom(M_1,M_2))\otimes M_1\cong \iHom(M_2,M_3)\otimes
(\iHom(M_1,M_2)\otimes M_1)\stackrel{id\otimes ev_{M_1,M_2}}{\longrightarrow}$$
$$\stackrel{id\otimes ev_{M_1,M_2}}{\longrightarrow}
\iHom(M_2,M_3)\otimes M_2\stackrel{ev_{M_2,M_3}}{\longrightarrow}M_3$$
which produces the {\em multipication morphism} 
$$\iHom(M_2,M_3)\otimes \iHom(M_1,M_2)\to \iHom(M_1,M_3).$$ 
 It is straightforward to check that this multiplication is
associative and compatible with the isomorphisms of Lemma \ref{can}.

Now let us fix an object $M\in \cM$. The multiplication morphism defines a
structure of an algebra on $A:=\iHom(M,M)$. Consider the category 
$\Mod_{\cC}(A)$ of right $A-$modules in the category $\cC$. The category
$\Mod_{\cC}(A)$ has an obvious structure of a left module category over $\cC$.
It is easy
to see that the functor $\iHom(M,?): \cM \to \Mod_{\cC}(A)$ has a natural
structure of module functor (this structure is induced by isomorphism (4) of
Lemma \ref{can}). We will say that $M\in \cM$ {\em generates} $\cM$ for any
$N\in \cM$ there is $X\in \cC$ such that $\Hom(X\otimes M,N)\ne 0$. It is easy
to see that in the case of exact module category $\cM$ the object $M$ 
generates $\cM$ if and only if its simple subquotients represent all 
equivalence classes in $Irr(\cM)_i$.

\begin{theorem} \label{morita}
Let $\cM$ be an exact module category over $\cC$ and assume
that $M\in \cM$ generates $\cM$. Then the functor $\iHom(M,?): \cM \to
\Mod_{\cC}(A)$ where $A=\iHom(M,M)$ is an equivalence of module categories.
\end{theorem}

\begin{proof} Note that in the case of an exact module category $\cM$ the 
functor $\iHom(M,?)$ is exact. The rest of the proof is parallel to the 
proof of Theorem 3.1 in \cite{O}. \end{proof}

\begin{definition} We will say that an algebra $A\in \cC$ is exact if the
category $\Mod_{\cC}(A)$ is exact.\end{definition}

\begin{example}\label{ihom} 
It is instructive to calculate $\iHom$ for the category
$\Mod_{\cC}(A)$. Let $M,N\in \Mod_{\cC}(A)$. Then $M^*$ has a natural
structure of a left $A^{**}-$module and $\iHom(M,N)=N\otimes^AM^*$ where
$N\otimes^AM^*:=(M\otimes_A{}^*N)^*$ (note that ${}^*N$ has a natural
structure of left $A-$module). Thus $N\otimes^AM^*$ is naturally a
subobject of $N\otimes M^*$ while $N\otimes_A{}^*M$ is a quotient of
$N\otimes {}^*M$. We leave to the reader to state and prove the associativity
properties of $\otimes^A$. One deduces from this description of $\iHom$ that
exactness of $A$ is equivalent to biexactness of $\otimes^A$ (and to
biexactness of $\otimes_A$).  
\end{example} 

\subsection{Dual category} In this subsection we show that there exists
a good notion of the dual category with respect to an exact module category.

Let $\cM_1$ and $\cM_2$ be two exact module categories over $\cC$.
Note that the category $Fun_{\cC}(\cM_1,\cM_2)$ of the additive module 
functors from $\cM_1$ to $\cM_2$ is abelian (note that such functors are
automatically of finite length).  

\begin{lemma} Let $\cM_1, \cM_2, \cM_3$ be exact module categories over 
$\cC$. The bifunctor of composition $Fun_{\cC}(\cM_2,\cM_3)\times Fun_{\cC}
(\cM_1,\cM_2)\to Fun_{\cC}(\cM_1,\cM_3)$ is biexact. \end{lemma}

\begin{proof} This is an immediate consequence of Proposition \ref{flat}.
\end{proof}

Another immediate consequence of Proposition \ref{flat} is the following:

\begin{lemma} Let $\cM_1, \cM_2$ be exact module categories over $\cC$. Any 
functor $F\in Fun_{\cC}(\cM_1,\cM_2)$ has both right and left adjoint.
\end{lemma}

Observe that an adjoint to a module functor has a natural structure of the 
module functor (we leave for the reader to define this). In particular,
it follows that the category $Fun_{\cC}(\cM,\cM)$
is a rigid monoidal category. We denote this category as $\cC^*_{\cM}$ and call
it the dual to $\cC$ with respect to $\cM$. 

We also have the following immediate

\begin{corollary} Let $\cM_1, \cM_2$ be exact module categories over $\cC$. 
Any functor $F\in Fun_{\cC}(\cM_1,\cM_2)$ maps projective objects to 
projectives. \end{corollary}

In view of Example \ref{exex} (ii) this Corollary is a generalization of 
Theorem \ref{projsurj} (but this does not give a new proof of Theorem 
\ref{projsurj}).

\begin{proposition} \label{finite}
The category $Fun_{\cC}(\cM_1,\cM_2)$ is finite. In
particular, the category $\cC^*_\cM$ is finite.
\end{proposition}

\begin{proof} We are going to use Theorem \ref{morita}. Thus $\cM_1=
\Mod_{\cC}(A_1)$ and $\cM_2=\Mod_{\cC}(A_2)$ for some algebras $A_1, A_2\in 
\cC$. It is easy to see that the category $Fun_{\cC}(\cM_1,\cM_2)$ is 
equivalent to the category of $(A_1,A_2)-$bimodules. But this category clearly
has enough projective objects: for any projective $P\in \cC$ the bimodule 
$A_1\otimes P\otimes A_2$ is projective. \end{proof}

\begin{lemma}\label{unit}
 The unit object $\be \in \cC^*_{\cM}$ is a direct sum of 
projectors to subcategories $\cM_i$. Each such projector is a simple object.
\end{lemma}

\begin{proof} The first statement is clear. For the second statement it is 
enough to consider the case when $\cM$ is indecomposable. Let $F$ be a nonzero
module 
subfunctor of the identity functor. Then $F(X)\ne 0$ for any $X\ne 0$. 
Hence $F(X)=X$ for any simple $X\in \cM$ and thus $F(X)=X$ 
for any $X\in \cM$ since $F$ is exact. \end{proof}

Thus the category $\cC^*_{\cM}$ is a finite multi-tensor category; in 
particular
if $\cM$ is indecomposable then $\cC^*_{\cM}$ is finite tensor category.
Note that by the definition $\cM$ is a module category over $\cC^*_{\cM}$. 

\begin{lemma} The module category $\cM$ over $\cC^*_{\cM}$ is exact.
\end{lemma}

\begin{proof} Let $A\in \cC$ be an algebra such that $\cM=Mod_{\cC}(A)$. Thus
the category $\cC^*_{\cM}$ is identified with the category $\Bimod(A)^{op}$ of
$A-$bimodules with opposite tensor product (because $A-$bimodules
act naturally on $\Mod_{\cC}(A)$ from the right).
Any projective object in the category of $A-$bimodules is a direct summand
of the object of the form $A\otimes P\otimes A$ for some projective $P\in \cC$.
Now for any $M\in Mod_{\cC}(A)$ one has that $M\otimes_AA\otimes P\otimes A=
(M\otimes P)\otimes A$ is projective by exactness of the category 
$Mod_{\cC}(A)$. The Lemma is proved.\end{proof}

\begin{example} \label{bihom}
It is instructive to consider the internal
Hom for the category $\Mod_{\cC}(A)$ considered as a module category over
$\cC^*_{\cM}=\Bimod(A)$. We leave to the reader to check that 
$\iHom_{\cC^*_{\cM}}(M,N)={}^*M\otimes N$ (the right hand side has obvious 
structure of $A-$bimodule). In particular $B=\iHom_{\cC^*_{\cM}}(A,A)=
{}^*A\otimes A$ is an algebra in the category of $A-$bimodules. Thus $B$ is
an algebra in the category $\cC$ and it is easy to see from definitions that
the algebra structure on $B={}^*A\otimes A$ comes from the evaluation morphism
$ev: A\otimes {}^*A\to \be$. Moreover, the coevaluation morphism induces
an imbedding of algebras $A\to {}^*A\otimes A\otimes A\to {}^*A\otimes A=B$ and
the $A-$bimodule structure of $B$ comes from the left and right multiplication
by $A$.
\end{example}

Thus for any exact module category $\cM$ over $\cC$ the category 
$(\cC^*_{\cM})^*_{\cM}$ is well defined. There is an obvious tensor functor
$can: \cC \to (\cC^*_{\cM})^*_{\cM}$.

\begin{theorem}\label{2star}
The functor $can: \cC \to (\cC^*_{\cM})^*_{\cM}$ is an
equivalence of categories.\end{theorem}

\begin{proof} Let $A$ be an algebra such that $\cM =\Mod_{\cC}(A)$. The
category $\cC^*_{\cM}$ is identified with the category $\Bimod(A)^{op}$.
The category $(\cC^*_{\cM})^*_{\cM}$ is identified with the category
of $B-$bimodules in the category of $A-$bimodules (here $B$ is the same
as in Example \ref{bihom} and is considered as an algebra in the category
of $A-$modules). But this latter category is tautologically identified with
the category of $B-$bimodules (here $B$ is an algebra in the category $\cC$)
since for any $B-$module one reconstructs the $A-$module structure via the
imbedding $A\to B$ from Example \ref{bihom}. We are going to use the
following

\begin{lemma} Any left $B-$module is of the form ${}^*A\otimes X$ for some
$X\in \cC$ with the obvious structure of an $A-$module. Similarly, any right 
$B-$module is of the form $X\otimes A$.\end{lemma}

\begin{proof} Let us consider $\cC$ as a module category over itself. 
Consider an object ${}^*A\in \cC$ as an object of this module category. Then
by Example \ref{ihom} $\iHom({}^*A,{}^*A)={}^*A\otimes A=B$ and the statement
follows from Theorem \ref{morita}. The case of right modules is completely
parallel. \end{proof}

It follows from the Lemma that any $B-$bimodule is of the form ${}^*A\otimes X
\otimes A$ and it is easy to see that $can(X)={}^*A\otimes X\otimes A$. The
Theorem is proved.
\end{proof}

\begin{corollary}\label{indecomp} Assume that $\cC$ is a finite tensor (not
only multi-tensor) category. Then an exact module category $\cM$ over $\cC$ is 
indecomposable over $\cC^*_\cM$.\end{corollary}

\begin{proof} This is an immediate consequence of Theorem \ref{2star} and
Lemma \ref{unit}. \end{proof}

Let $\cM$ be a fixed module category over $\cC$. For any other module category
$\cM_1$ over $\cC$ the category $\Fun_{\cC}(\cM_1,\cM)$ has obvious structure
of a module category over $\cC^*_{\cM}=\Fun_{\cC}(\cM,\cM)$. 

\begin{lemma} The module category $\Fun_{\cC}(\cM_1,\cM)$ over $\cC^*_{\cM}$ 
is exact. \end{lemma} 

\begin{proof} Assume that $\cM=\Mod_{\cC}(A)$ and $\cM_1=\Mod_{\cC}(A_1)$.
Identify $\cC^*_{\cM}$ with the category of $A-$bimodules and $\Fun_{\cC}
(\cM_1,\cM)$ with the category of $(A_1-A)-$bimodules. Any projective object
of $\Bimod(A)$ is a direct summand of an object of the form $A\otimes P\otimes
A$ for some projective $P\in \cC$. Let $M$ be an $(A_1-A)-$bimodule, then
$M\otimes_AA\otimes P\otimes A=M\otimes P\otimes A$. Now $\Hom_{A_1-A}
(M\otimes P\otimes A,?)=\Hom_{A_1}(M\otimes P,?)$ (here $\Hom_{A_1-A}$ is
the Hom in the category of $(A_1-A)-$bimodules and $\Hom_{A_1}$ is the Hom
in the category of left $A_1-$modules) and it is enough to check that
$M\otimes P$ is a projective left $A_1-$module. This is equivalent to
$(M\otimes P)^*$ being injective (since $N\mapsto N^*$ is an equivalence of
the category of left $A-$modules to the category of right $A-$modules).
But $(M\otimes P)^*=P^*\otimes M^*$ and results follows from projectivity
of $P^*$ and Lemma \ref{proinj}. 
\end{proof}

The proof of the following Theorem is similar to the proof of Theorem 
\ref{2star} and is left to the reader.

\begin{theorem} Let $\cM$ be an exact module category over $\cC$.
The maps $\cM_1\mapsto \Fun_{\cC}(\cM_1,\cM)$ and $\cM_2
\mapsto \Fun_{\cC^*_{\cM}}(\cM_2,\cM)$ are mutually inverse bijections
of the sets of equivalence classes of exact module categories over $\cC$ and 
over $\cC^*_{\cM}$. 
\end{theorem}

Following \cite{Mu1} we will say that the categories $\cC$ and 
$(\cC^*_{\cM})^{op}$ are {\em weakly Morita equivalent}. 

Let $\cM$ be an exact module category over $\cC$. For $X,Y\in \cM$ we have
two notions of internal Hom --- with values in $\cC$ and with values in 
$\cC^*_{\cM}$, denoted by $\iHom_\cC$ and $\iHom_{\cC^*_{\cM}}$ respectively.
The following simple consequence of calculations in Examples \ref{ihom} and
\ref{bihom} is very useful.

\begin{proposition}\label{basic}
 {\em (``Basic identity'')} Let $X,Y,Z \in \cM$. There is
a canonical isomorphism
{\em $$\iHom_\cC(X,Y)\otimes Z\simeq {}^*\iHom_{\cC^*_{\cM}}(Z,X)\otimes Y.$$}
\end{proposition}

\begin{proof} By Theorem \ref{2star} it is enough to find a canonical
isomorphism
$${}^*\iHom_\cC(Z,X)\otimes Y\simeq \iHom_{\cC^*_{\cM}}(X,Y)\otimes Z.$$
This isomorphism is constructed as follows.
Choose an algebra $A$ such that $\cM =\Mod_\cC(A)$. By Example
\ref{ihom} the LHS is ${}^*(X\otimes^AZ^*)\otimes Y={}^*(Z\otimes_A{}^*X)^*
\otimes Y=(Z\otimes_A{}^*X)\otimes Y$. On the other hand by Example \ref{bihom}
the RHS is $Z\otimes_A({}^*X\otimes Y)$. Thus the associativity isomorphism
gives a canonical isomorphism of the LHS and RHS. Observe that the isomorphism
inverse to the one we constructed is the image of the identity under the 
homomorphism
$$\Hom(Y,Y)\to \Hom(\iHom_{\cC^*_{\cM}}(X,Y)\otimes X, Y)\to$$
$$\Hom(\iHom_{\cC^*_{\cM}}(X,Y)\otimes \iHom_\cC(Z,X)\otimes Z, Y)\simeq$$
$$\Hom(\iHom_\cC(Z,X)\otimes \iHom_{\cC^*_{\cM}}(X,Y)\otimes Z, Y)\simeq$$
$$\Hom(\iHom_{\cC^*_{\cM}}(X,Y)\otimes Z, {}^*\iHom_\cC(Z,X)\otimes Y)$$
and thus does not depend on the choice of $A$. 
\end{proof}

\begin{remark} Following \cite{Mu1} one can construct from $\cM$ a 2-category
with 2 objects $\fA, \fB$ such that $\End(\fA)\cong \cC$, $\End(\fB)\cong
(\cC^*_{\cM})^{op}$, $\Hom(\fA, \fB)\cong \cM$, and $\Hom(\fB, \fA)=
Fun_{\cC}(\cM,\cC)$. In this language Proposition \ref{basic} expresses the 
associativity of the composition of Hom's. \end{remark}

\subsection{The Drinfeld double and weak Morita invariance of the 
Frobenius-Perron dimension} Let $\cC$ be a finite multi-tensor category and
let $\cM$ be a module category over $\cC$. Consider $\cM$ as a module
category over $\cC \boxtimes \cC^*_{\cM}$. Clearly this module category is
exact. Let $Z(\cC)$ denote the Drinfeld
center of the category $\cC$.

\begin{theorem}\label{double}
The category $(\cC \boxtimes \cC^*_\cM)^*_\cM$ is canonically
equivalent to $Z(\cC)$. In particular the category $Z(\cC)$ is finite.
\end{theorem}

\begin{proof} (see \cite{Ob}) Any object of $(\cC \boxtimes \cC^*_\cM)^*_\cM$
commutes with $\cC^*_\cM-$action, so is an object $X$ of $\cC$ by Theorem 
\ref{2star}. Additionally any object of $(\cC \boxtimes \cC^*_\cM)^*_\cM$
commutes with $\cC$ whence we get a structure of the object of $Z(\cC)$ on
$X$. We leave to the reader to check that commutative diagrams from 
definitions of $(\cC \boxtimes \cC^*_\cM)^*_\cM$ and $Z(\cC)$ correspond. 

The second assertion is a special case of Proposition \ref{finite}.
\end{proof}

\begin{corollary}\label{morZ} There is a canonical equivalence 
$Z(\cC)\simeq Z(\cC^*_\cM)$.\end{corollary}

\begin{proof} Theorem \ref{double} is symmetric in $\cC$ and $\cC^*_\cM$. Thus
both $Z(\cC)$ and $Z(\cC^*_\cM)$ are canonically equivalent to 
$(\cC \boxtimes \cC^*_\cM)^*_\cM$. \end{proof} 

\begin{remark} It is easy to see that under the equivalence of Corollary
\ref{morZ} the braiding of the category $Z(\cC)$ corresponds to the
inverse braiding of the category $Z(\cC^*_{\cM})$. In other words we have
an equivalence of braided categories $Z(\cC)\simeq Z((\cC^*_{\cM})^{op})$.
\end{remark} 

We are going to use the particular case of Theorem \ref{double} with $\cM=\cC$.
In this case $\cC^*_\cM=\cC^{op}$ and we get

\begin{corollary} There is a canonical equivalence $Z(\cC)\simeq (\cC \boxtimes
\cC^{op})^*_\cC$. \end{corollary}

Let $F: Z(\cC)\to \cC$ denote the canonical forgetful functor and let
$I: \cC \to Z(\cC)$ denote the right adjoint functor of $F$ (thus
$\Hom(F(X),Y)=\Hom(X,I(Y))$). The functor $I$ can be expressed in terms of
the internal Hom in the following way:

\begin{lemma}\label{Iihom}
 We have canonically $I(X)=\iHom_{Z(\cC)}(\be,X)$.\end{lemma}

\begin{proof} By definition $F(X)=F(X)\otimes \be =X\otimes \be$ (in the last
equation $\be$ is an object of the module category $\cC$ over $Z(\cC)$).
Thus $\Hom(F(X),?)=\Hom(X,\iHom_{Z(\cC)}(\be,?))$ and the Lemma is proved.
\end{proof}

\begin{proposition}\label{facts} {\em (i)} The functor $F$ is surjective. 

{\em (ii)} The functor $I$ is exact.

{\em (iii)} For $X\in \cC$ and $Y\in Z(\cC)$ we have canonical isomorphisms
$$I(X\otimes F(Y))\simeq I(X)\otimes Y\;\;\; \mbox{and}\;\;\; I(F(Y)\otimes X)
\simeq Y\otimes I(X).$$ \end{proposition}

\begin{proof} (i) is an immediate consequence of Corollary \ref{indecomp};
(ii) follows from Lemma \ref{Iihom} and Corollary \ref{biex}; (iii) is a
special case of Lemma \ref{can}. \end{proof}

Recall that in Section \ref{FPdim} the Frobenius-Perron dimensions were
defined. 

\begin{lemma}\label{frac} For any object $V\in \cC$ one has $d_+(I(V))=
\frac{d_+(Z(\cC))}{d_+(\cC)}d_+(V)$.\end{lemma}

\begin{proof} Recall the virtual object $R_\cC$ from Section \ref{FPdim}. It
follows from Proposition \ref{facts} (i) and formula $(1)$ that 
$F(R_{Z(\cC)})=\frac{d_+(Z(\cC))}
{d_+(\cC)}R_\cC$. Note that $\dim \Hom(R_\cC,X)$ is well defined for any $X$
since $R_\cC$ is ``virtually projective''. One has
\begin{eqnarray*}
d_+(I(V))
&=& \oplus_{X\in Irr(Z(\mC))}d_+(X)[I(V):X]\\
&=& \oplus_{X\in Irr(Z(\mC))}d_+(X)\dim \Hom(P(X),I(V)) \\
&=& \dim \Hom(R_{Z(\cC)},I(V))=\frac{d_+(Z(\cC))}{d_+(\cC)}
\dim \Hom(R_\cC,V)\\ &=&\frac{d_+(Z(\cC))}{d_+(\cC)}d_+(V).
\end{eqnarray*}
The Lemma is proved. \end{proof}

\begin{lemma}\label{one} We have $d_+(I(\be))=d_+(\cC)$. \end{lemma}

\begin{proof} We have by Lemma \ref{Iihom}:
$$d_+(I(\be))=d_+(F(I(\be))=d_+(\iHom_{Z(\cC)}(\be,\be)
\otimes \be).$$
Now by Proposition \ref{basic} $\iHom_{Z(\cC)}(\be,\be)\otimes \be =
{}^*\iHom_{\cC 
\boxtimes \cC^{op}}(\be,\be)\otimes \be$ and thus $d_+(I(\be))=
d_+(\iHom_{\cC \boxtimes \cC^{op}}(\be,\be))$. The simple objects
of $\cC \boxtimes \cC^{op}$ are of the form $X\boxtimes Y$ where $X,Y\in 
Irr(\cC)$ and their projective covers are of the form $P(X)\boxtimes P(Y)$.
Hence
$$\begin{array}{cl}
&d_+(\iHom_{\cC \boxtimes \cC^{op}}(\be,\be))\\=&
\sum_{X,Y\in Irr(\cC)}
d_+(X)d_+(Y)[\iHom_{\cC \boxtimes \cC^{op}}(\be,\be):
X\boxtimes Y]\\
=&\sum_{X,Y\in Irr(\cC)}d_+(X)d_+(Y)\dim \Hom
(P(X)\boxtimes P(Y),\iHom_{\cC \boxtimes \cC^{op}}(\be,\be))\\ =&
\sum_{X,Y\in Irr(\cC)}d_+(X)d_+(Y)\dim \Hom(P(X)\otimes P(Y),
\be)\\ =&
\sum_{X,Y\in Irr(\cC)}d_+(X)d_+(Y)\dim \Hom(P(X),P(Y)^*)\\ =&
\sum_{X,Y\in Irr(\cC)}d_+(X)d_+(Y)[P(Y)^*:X]\\ =&
\sum_{Y\in Irr(\cC)}d_+(Y)d_+(P(Y)^*)=d_+(\cC).
\end{array}$$
The Lemma is proved. \end{proof}

\begin{theorem}\label{square} We have $d_+(Z(\cC))=(d_+(\cC))^2$. \end{theorem}

\begin{proof} We have $d_+(I(\be))=\frac{d_+(Z(\cC))}{d_+(\cC)}$ by
Lemma \ref{frac} and $d_+(I(\be))=d_+(\cC)$ by Lemma \ref{one}.
\end{proof}

\begin{corollary} For any exact module category $\cM$ we have
$d_+(\cC)=d_+(\cC^*_\cM)$. \end{corollary}

\begin{proof} By Corollary \ref{morZ} and Theorem \ref{square} we have
$(d_+(\cC))^2=(d_+(\cC^*_\cM))^2$. Since both numbers $d_+(\cC)$ and
$d_+(\cC^*_\cM)$ are positive we are done. \end{proof}

\subsection{Dualization of tensor functors}

Let $\cC,\cD$ be finite multi-tensor categories, 
$\cM$ be an exact module category over $\cD$, and 
$F:\cC\to \cD$ be a tensor functor
(i.e., we require that $F(\be)=\be$). Then $\cM$ is 
a module category over $\cC$ which is, obviously, not always
exact (e.g., $\cC$ is trivial, $\cM=\cD$). 

\begin{definition}
The pair $(F,\cM)$ is called an exact pair if 
$\cM$ is exact over $\cC$.  
\end{definition}

Suppose $(F,\cM)$ is an exact pair. 
Then we have the obvious dual tensor functor $F^*: \cD_\cM^*\to
\cC_\cM^*$, and $(F^*,\cM)$ is an exact pair. We say that 
the exact pair $(F^*,\cM)$ is dual to the 
pair $(F,\cM)$, and write $(F^*,\cM)=(F,\cM)^*$. 
Clearly, for any exact pair $T$, one has $T^{**}=T$. 

For simplicity, from now till the end of the subsection 
we will consider only exact pairs $T$ in which $\cC$, $\cD$
are tensor (i.e. not just multi-tensor) categories, and 
$\cM$ is indecomposable over $\cC$
(the class of such pairs are obviously stable under 
dualization). We note, however, that 
the results below can be extended to the general case. 

\begin{definition}
An exact pair is surjective if $F$ is surjective,
and injective if $F$ is injective.
\end{definition}

\begin{theorem}\label{dualiz}
The dualization map takes surjective exact pairs into injective 
ones, and vice versa. 
\end{theorem}

\begin{proof} Let $T=(F:\cC\to \cD,\mM)$ be an exact pair. 
Let $d_+(\cC)=c$, $d_+(\cD)=d$.
Then $d_+(\cC_\cM^*)=c$, $d_+(\cD_\cM^*)=d$.
 
Assume first that $F$ is
injective, but $F^*$ is not surjective. 
Then by Proposition \ref{injsurj}, 
$d_+({\rm Im}F^*)<c$. Since $F$ 
is injective, we also have $d_+({\rm Im}F^*)<d$ (as $c\le d$).
The functor $F$ factors through ${\mathcal E}=({\rm
Im}F^*)^*_{\cM}$ (it is not difficult to show that 
$\cM$ is exact and indecomposable over 
${\rm Im}F^*$, so ${\mathcal E}$ is a finite tensor category). 
Since $d_+({\mathcal E})=d_+({\rm Im}F^*)<{\rm min}(c,d)$, 
by Proposition \ref{notinjsurj}, $F$ is not
injective. Contradiction. 

Assume now that $F$ is
surjective, but $F^*$ is not injective. 
Then by Proposition \ref{injsurj}, 
$d_+({\rm Im}F^*)<d$. Since $F$ 
is surjective, we also have $d_+({\rm Im}F^*)<c$.
The functor $F$ factors through ${\mathcal E}=({\rm
Im}F^*)^*_{\cM}$. 
Since $d_+({\mathcal E})=d_+({\rm Im}F^*)<{\rm min}(c,d)$, 
by Proposition \ref{notinjsurj}, $F$ is not
surjective. Contradiction. 
\end{proof}

\subsection{Lagrange's theorem for finite tensor categories}

\begin{theorem}\label{Lagr} 
Let $\mD$ be a finite tensor category, and
$\mC\subset \mD$ be a tensor subcategory. 
Then $d_+(\mD)/d_+(\mC)$ is an algebraic integer.
\end{theorem}

\begin{proof} Consider the natural embedding $F: \mC\boxtimes
\mD^{op}\to \mD\boxtimes \mD^{op}$.  
Consider $\mM=\mD$ as a module category over $\mD\boxtimes
\mD^{op}$. It is easy to check that the pair $(F,\mM)$ is exact,
and $\mM$ is indecomposable over $\mC\boxtimes \mD^{op}$.
Thus, Theorem \ref{dualiz} applies, and the functor 
$F^*: (\mD\otimes \mD^{op})_\mM^*=Z(\mD)\to
(\mC\otimes \mD^{op})_\mM^*$ is surjective. The dimension of the
first category is $d_+(\mD)^2$ and of the second 
one $d_+(\mC)d_+(\mD)$, so as explained in \ref{freen}, 
$d_+(\mD)/d_+(\mC)$ is an algebraic integer. We are done. 
\end{proof}

\begin{corollary} \cite{Sch2}
The dimension of any quasi-Hopf quotient of a quasi-Hopf algebra 
divides the dimension of this quasi-Hopf algebra. 
\end{corollary}

\section{Examples} 
In this section we present some cases when we were able to classify exact
module categories. Let $k$ denote an algebraically closed field.

\subsection{Finite groups} Let $G$ be a finite group. Consider the tensor
category $\Rep_k(G)$ of representations of $G$ over the field $k$. Let
$H\subset G$ be a subgroup and $\psi \in H^2(H,k^*)$. A choice of a
cocycle representing $\psi$ defines a central extension
$$1\to k^*\to \tilde H\to H\to 1.$$
Let $\Rep_k(H,\psi)$ denote the category of representations $\tilde H$
over $k$ such that $z\in k^*$ acts via multiplication by $z$. Clearly
$\Rep_k(H,\psi)$ is a module category over $\Rep_k(G)$. Obviously, the
module category $\Rep_k(H,\psi)$ is exact and it is easy to see that it
does not depend on a choice of cocycle representing $\psi$. 

\begin{proposition}\label{repg}
The indecomposable exact module categories over $Rep_k(G)$
are of the form $\Rep_k(H,\psi)$ and are classified by conjugacy classes of 
pairs $(H,\psi)$ where $H\subset G$ is a subgroup and $\psi \in H^2(H,k^*)$.
\end{proposition}

\begin{proof} We begin with the following general

\begin{lemma}\label{rightsimple}
 Let $\cM$ be an exact module category over a finite tensor
category $\cC$ and let $M\in \cM$ be a simple object. The algebra $A=
\iHom(M,M)$ in the category $\cC$ has no nontrivial right ideals (``simple
from the right''). \end{lemma}

\begin{proof} The module category $\cM$ is module equivalent to the category
$\Mod_\cC (A)$ and under this equivalence the object $M$ corresponds to $A$
considered as a right $A-$module. Since $M$ is simple the result follows.
\end{proof} 

For example let $V\in \Rep_k(H,\psi)$ be a simple object. One calculates 
immediately that $\iHom(V,V)=\Ind_H^G(\End(V))$ and thus the algebra
$A=\Ind_H^G(\End(V))$ is simple from the right.

\begin{lemma}\label{Galg}
Let $A\in \Rep_k(G)$ be a simple from the right algebra. Then
$A$ is of the form $A=\Ind_H^G(\End(V))$ where $V$ is a simple object of
$\Rep_k(H,\psi)$ for some $H\subset G$ and $\psi \in H^2(H,k^*)$.
\end{lemma}

\begin{proof}  
In other words $A$ is an associative algebra with an action of $G$ by 
automorphisms and without nontrivial $G-$invariant right ideals. Since a 
group action preserves the radical, the algebra $A$ is semisimple. The group 
$G$ acts transitively on the set of minimal central idempotents of $A$ (since 
otherwise we would have a $G-$invariant direct summand). Let $H$ be the 
stabilizer of a minimal central idempotent $e$; clearly $A=\Ind_H^G(eAe)$ and 
$eAe$ is a matrix algebra. Thus $eAe=\End(V)$ where $V$ is a projective 
representation of $H$; the representation $V$ is irreducible since otherwise 
the $G-$span of the annihilator of an $H-$submodule in $V$ would be a 
$G-$invariant right ideal. \end{proof}

Since for any object $M\in \cM$ the algebra $A=\iHom(M,M)$ determines the
exact module category $\cM$ uniquely, the proposition is proved.
\end{proof}

\begin{remark} (i) Proposition \ref{repg} is new only in the case $\chara (k)
>0$, see e.g. \cite{O}. Actually, our proof repeats the characteristic 0 proof.

(ii) There is another proof of Proposition \ref{repg} along the lines of
\cite{Ob}. One can easily show that many other results of \cite{Ob} remain
true in positive characteristic in the setting of exact module categories,
for example the classification of module categories over the Drinfeld double
of a finite group.

(iii) It seems plausible that the converse to Lemma \ref{rightsimple} is
true, that is for a simple from the right algebra $A\in \cC$ the module
category $\Mod_{\cC}(A)$ is exact. \end{remark} 

\subsection{Finite supergroups} In this section we will assume that $\chara(k)
\ne 2$. Let $G$ be a finite group, $W$ be a 
representation (possibly zero) of $G$ and $u\in G$ be a central element of 
order $\le 2$ acting by $(-1)$ on $W$. Regard $W$ as an odd supervector space
and consider the supergroup $G\ltimes W$. Let us consider the category
$\Rep (G\ltimes W,u)$ of representations on super vector spaces $V$ of 
$G\ltimes W$ such that $u$
acts on $V$ via the parity automorphism. The category $\Rep (G\ltimes W,u)$ has
an obvious structure of a tensor category. Recall that according to P.~Deligne 
\cite{De} in a case $\chara(k)=0$ the most general finite {\em symmetric} 
tensor category is of the form $\Rep (G\ltimes W,u)$. 

One can construct exact module categories over $\Rep (G\ltimes W,u)$ in the
following way. Let $H\subset G$ be a subgroup and let $Y$ be an $H-$invariant
subspace of $W$. Let $B$ be an $H-$invariant quadratic form on $Y$ (possibly 
degenerate) and let $\Cl(Y,B)$ denote the corresponding Clifford algebra. Let 
$\psi \in Z^2(H,k^*)$ be a two cocycle and let $k[H]_\psi$ denote the
corresponding twisted group algebra. Let $\cM^0(Y,B,H,\psi)$ denote the
category of finite dimensional super vector spaces $Z$ endowed 
with the following structures: 

(i) $Z$ is a $\Cl(Y,B)-$module, that is for any $v\in Y$ we have an odd 
endomorphism $\p_v$ of $Z$ such that $\p_v^2=B(v,v)$;

(ii) $Z$ is $k[H]_\psi-$module, that is for any $h\in H$ we have an 
automorphism $h$ of $Z$ and $h_1\cdot h_2=\psi(h_1,h_2)h_1h_2$;

(iii) Structures (i) and (ii) are compatible: we have $h\p_vz=\p_{hv}hz$ for
all $h\in H, v\in Y, z\in Z$.

In other words $\cM^0(Y,B,H,\psi)$ is the category of (super) representations 
of the suitably defined smash-product $k[H]_\psi \ltimes \Cl(Y,B)$.

The category $\cM^0(Y,B,H,\psi)$ has a natural structure of a module category
over $\Rep (G\ltimes W,u)$: for $S\in \Rep (G\ltimes W,u)$ and $Z\in
\cM^0(Y,B,H,\psi)$ we set $S\otimes Z$ to be the usual tensor product of vector
spaces with the following action of $\Cl (Y,B)$ and $k[H]_\psi$:
$$\p_v(s\otimes z)=vs\otimes z+(-1)^{|s|}s\otimes \p_v(z);
\;\; h(s\otimes z)=h(s)\otimes h(z).$$

We leave for the reader to check that the module category $\cM^0(Y,B,H,\psi)$
over $\Rep (G\ltimes W,u)$ is exact (note that being projective in the
category $\cM^0(Y,B,H,\psi)$ is the same as being free over $\wedge (\Ker(B))
\subset \Cl(Y,B)$ and being projective over $k[H]_\psi$). Note that the
module category $\cM^0(Y,B,H,\psi)$ is not always indecomposable: in a case
when $\dim (Y/\Ker(B))$ is even it is equivalent to the direct sum of two
indecomposable module categories both being equivalent to the category
defined similarly to $\cM^0(Y,B,H,\psi)$ but with $Z$ assumed to be a 
usual (not super) vector space. Thus we set
$$\cM(Y,B,H,\psi)=\left\{ \begin{array}{ll}\cM^0(Y,B,H,\psi)&\mbox{if}\; 
\dim (Y/\Ker(B))\; \mbox{is odd};\\ 
\mbox{an indecomposable summand}&\mbox{if}\; 
\dim (Y/\Ker(B))\; \mbox{is even}.\\
\mbox{of}\; \cM^0(Y,B,H,\psi)&\end{array}\right.$$

Note that in the case $\dim (Y/\Ker(B))$ is odd the Clifford algebra
$\Cl(Y,B)$ has a unique irreducible super representation $S$ and in the
case $\dim (Y/\Ker(B))$ is even the Clifford algebra $\Cl(Y,B)$ has a
unique irreducible reprsentation $S$. Clearly in both cases the group $H$ 
acts on $S$ projectively; let us choose a corresponding 2-cocycle $\psi_0$
(in other words: the group $H$ maps to the orthogonal group $O(Y/\Ker(B))$
and $\psi_0$ is the inverse image of a cocycle defining the spinor group).
It is easy to see that the simple objects of $\cM(Y,B,H,\psi)$ are of the
form $V\otimes S$ where $V$ is an irreducible projective representation
of $H$ corresponding to the 2-cocycle $\psi-\psi_0$ (here $\Cl(Y,B)$ acts
trivially on the first factor and $H$ acts diagonally). The quadratic form
$B$ induces a non-degenerate quadratic form on $Y/\Ker(B)$, hence a
non-degenerate quadratic form on $(Y/\Ker(B))^*$ and consequently a
quadratic form (denoted by the same letter $B$) on $(W/\Ker(B))^*$. The
corresponding Clifford algebra $\Cl((W/\Ker(B))^*,B)$ has an obvious
action of $\hat H\ltimes W$ where $\hat H\subset G$ denote the subgroup 
generated by $H$ and $u$. One calculates readily
$$\iHom_{\cM(Y,B,H,\psi)}(V\otimes S,V\otimes S)=\Ind_{\hat H}^G(
\Ind_H^{\hat H}(\End(V))\otimes \Cl((W/\Ker(B))^*,B))$$ 
(tensor product here is in the super sense).

The main result of this subsection is the following

\begin{theorem}\label{super}
Any indecomposable exact module category over $\Rep 
(G\ltimes W,u)$ is of the form $\cM(Y,B,H,\psi)$. Two module categories
$\cM(Y,B,H,\psi)$ and $\cM(Y',B',H',\psi')$ are module equivalent if
and only if there is $g\in G$ such that $g(Y)=Y'$, $g(B)=B'$, $gHg^{-1}=H'$,
and $g(\psi)\psi'^{-1}$ is a coboundary.
\end{theorem}

\begin{proof} We are going to proceed in the same way as in the proof of
Proposition \ref{repg}.

Thus we are going to classify simple from the right 
algebras in the category
$\Rep (G\ltimes W,u)$. In down to earth terms we are looking for finite
dimensional associative algebras $A$ (with unit) with the following structures:

(i) The group $G$ acts by automorphisms on $A$; in particular the element $u$
determines the structure of a superalgebra on $A$.

(ii) For any vector $v\in W$ we have an odd derivation $\p_v$ of $A$; the
assignment $v\mapsto \p_v$ is linear; for $v,w\in W$ we have $\p_v\p_w=-\p_w
\p_v$; in particular $(\p_v)^2=0$.

(iii) The structures (i) and (ii) are compatible; that is for any $g\in G, v\in
W, a\in A$ we have $g\p_v(a)=\p_{gv}(ga)$.

(iv) The algebra $A$ has no nontrivial right ideals $I$ such that 
$G(I)\subset I$ and $\p_v(I)\subset I$ for all $v\in W$.

Define inductively a filtration on $A$: set $A_{-1}=0$ and $A_i=\{ a\in A|
\p_v(a)\in A_{i-1} \forall v\in W\}$. Clearly $A_i\subset A_{i+1}$. Observe 
that

(a) For any $v\in W$ we have $\p_v(A_i)\subset A_{i-1}$;

(b) The filtration is $G-$invariant $G(A_i)\subset A_i$;

(c) The filtration is multiplicative $A_iA_j\subset A_{i+j}$;

(d) The filtration is exhausting $A_{\dim W}=A$.

(e) For any $a\in A_i\setminus A_{i-1}$ there exists $v\in W$ such that 
$\p_v(a)\in A_{i-1}\setminus A_{i-2}$.

In particular $A_0$ is a $G-$invariant subalgebra of $A$. Note that for any
nonzero right $G-$invariant ideal $I_0\subset A_0$ the right ideal $I_0A$ of 
$A$ is nonzero $G-$invariant and $\p_v-$invariant for any $v\in W$. Thus by
condition (iv) we have $I_0A=A$. In particular this applies to the radical
$R_0$ of $A_0$ (note that $R_0$ is automatically $G-$invariant) and hence
$R_0A=A$. But since $R_0$ is a nilpotent ideal, this equality is impossible and
thus $R_0=0$. Hence the algebra $A_0$ is semisimple and the filtration
$A_0\subset A_1\subset \ldots$ splits as a filtration of $A_0-$modules.
Hence $I_0A=A$ implies that $I_0=I_0A_0=A_0$ for any right ideal $I_0\subset
A_0$. Summarizing we get

(f) The algebra $A_0$ has no nontrivial $G-$invariant ideals. 

Hence there is a subgroup $H\subset G$ and an irreducible projective
representation $V$ of $H$ such that $A_0=\Ind_H^G(\End(V))$. 
Let $e_1, \ldots e_n$ be the even minimal central idempotents of the algebra 
$A_0$ (so $n=|G/\hat H|$ where $\hat H\subset G$ is the subgroup of $G$
generated by $H$ and $u$; thus $|\hat H/H|$ is 1 or 2). Using property (e)
one proves by an easy induction that $A=\bigoplus_{i=1}^ne_iAe_i$. Thus
$e_i$ are the central idempotents in the algebra $A$ and $A=
\Ind_{\hat H}^G(e_1Ae_1)$.
Thus we can (and will) assume that either $G=H$ or $G=H\times \{ 1,u\}$.

Let $W'=\{ v\in W | \p_v|_{A_1}=0\}$. It is clear from (b) that $W'$ is a
$G-$invariant subspace of $W$. It follows from (e) and (ii) that for 
$v\in W'$ we have $\p_v(a)=0$ for any $a\in A$. Let us denote $W/W'=X$.
For any $v\in X$ the derivation $\p_v$ of the algebra $A$ is well defined. 
Now let us define $U=\{ x\in A_1|\; x\; \mbox{is odd}, xa=u(a)x\; \forall 
a\in A_0\}$. 

\begin{lemma} The multiplication induces an isomorphisms $U\otimes A_0\to 
A_1/A_0$ and $A_0\otimes U\to A_1/A_0$.\end{lemma}

\begin{proof} Let us consider $A_0$ as an algebra with an action of the 
element $u$. Then $A_1/A_0$ is an $A_0-$bimodule with an action of $u$, and 
the two structures are compatible in the obviuos sense. It is clear that any 
such bimodule is a direct sum of simple bimodules, and it is enough to check 
the statement of the Lemma only for simple summands of $A_1/A_0$. 

Case 1: $G=H$. In this case $A_0$ is just the matrix algebra 
$A_0=\End(T_+\oplus T_-)$ where $u|_{T_+}=1$ and $u|_{T_-}=-1$. There are 
two simple $A_0-$bimodules with a $u-$action: $A_0$ itself and $A_0$ with 
the opposite parity. In both cases the space of $x$ such that $xa=u(a)x$ for 
all $a\in A_0$ is one dimensional and generated by $x_0=1|_{T_+}\oplus 
(-1)|_{T_-}$. On the other hand there is $v\in X$ such that $\p_v(x)\ne 0$. 
Then $\p_v(x_0)a=a\p_v(x_0)$ for all $a\in A_0$ and thus $\p_v(x_0)$ is 
proportional to $1\in A_0$ and therefore is even. So we can assume that $x_0$ 
is odd. Since $x_0$ is invertible, the Lemma follows in this case.

Case 2: $G=H\times \{ 1,u\}$. In this case $A_0$ is a sum of two matrix 
algebras permuted by $u$. There is only one simple $A_0-$bimodule with 
a $u-$action occuring in $A_1/A_0$, and it is easy to check the Lemma in
this case.
\end{proof}    

Note that for any $x\in U$ and $v\in X$ one has $\p_v(x)a=a\p_v(x)$ for any
$a\in A_0$; also $\p_v(x)$ is even. Thus $\p_v(x)$ is proportional to 
$1\in A_0$. Thus we have a well defined pairing $\beta : X\times U\to k$
such that $\p_v(x)=\beta(v,x)1$, and it is easy to see from definitions that
this pairing is non-degenerate. Thus we can (and will) identify the space
$U$ with the space $X^*$. Any element $a\in A_i$ determines a map 
$\wedge^i(X)\to A_0$ defined by $v_1\wedge \ldots \wedge v_i\mapsto \p_{v_1}
\ldots \p_{v_i}(a)$ or, equivalently, an element of $(\wedge^i(X))^*
\otimes A_0$. This element is 0 if and only if $a\in A_{i-1}$. On the other 
hand, multiplying elements of $U$ and $A_0$, one can construct an element 
$a\in A_i$ producing any given element of $(\wedge^i(X))^*\otimes A_0$. Thus
we have

(g) The algebra $A$ is generated by $A_1$.

Now let $x\in U$. For any $v\in X$ one has $\p_v(x^2)=\p_v(x)x-x\p_v(x)=0$.
Thus $x^2\in A_0$. Moreover, for any $a\in A_0$ we have $ax^2=x^2a$ and since
$x^2$ is even we see that $x^2$ is proportional to $1\in A_0$. Thus there
exists a quadratic form $B$ on $U=X^*$ such that $x^2=B(x,x)1$ for any 
$x\in U$. Clearly, the form $B$ is $G-$invariant. Thus $A$ is a quotient of 
$A_0\otimes \Cl(X^*,B)$ where $\Cl(X^*,B)$ is the Clifford algebra constructed
from the vector space of generators $X^*$ and the quadratic form $B$ (here
$\otimes$ is a tensor product of superalgebras). On the other hand it is easy
to see that the algebra $A_0\otimes \Cl(X^*,B)$ has no nontrivial ideals 
invariant under $G\ltimes W$.

Therefore we have identified an arbitrary simple from the right algebra in the 
category
$\Rep(G\ltimes W,u)$ with the internal Hom algebra of some simple object
in the category $\cM(Y,B,H,\psi)$. It follows from the above that the
algebras we constructed are pairwise nonisomorphic, whence we deduce the
second statement of the Theorem. The Theorem is proved.
\end{proof}

\begin{example} Consider the case when $G=\BZ/2\BZ$ and $u\in G$ is the 
nontrivial
element. In this case $\Rep(G\ltimes W,u)$ is just the category $\Rep(W)$ of 
representations of the supergroup $W$. There are two kinds of indecomposable
module categories over $\Rep(W)$: with two simple objects and with one simple
object (note that the category $\Rep(W)$ has just two simple objects and
both of them are invertible). The exact module categories of both kinds are
classified by a subspace $Y\subset W$ and a quadratic form $B\in S^2(Y^*)$;
the module category is semisimple if and only if the form $B$ is nondegenerate.
\end{example}

Note that from Theorem \ref{super}, one can obtain the
classification of fiber functors on $\Rep(G\ltimes W,u)$, i.e.
module categories which are equivalent to the category of vector
spaces. Namely, it is easy to see that the category 
$\cM(Y,B,H,\psi)$ is equivalent to the category of vector spaces 
if and only if the cocycle $\psi$ is nondegenerate (i.e. 
the algebra $k[H]_\psi$ is simple), and $B$ is a nondegenerate
quadratic form. In this case, the cocycle $\psi$ defines 
an irreducible projective representation $V$ of $H$. 

Recall now that the category $\mC:=\Rep(G\ltimes W,u)$ is symmetric, 
and that it admits a unique fiber functor which preserves the
symmetric structure (\cite{De}). This implies that 
equivalence classes of fiber functors 
on $\mC$ are in bijection with isomorphism classes
of triangular Hopf algebras $A$ such that $\Rep A=\mC$. 
Thus, Theorem \ref{super} implies that triangular Hopf algebras
$A$ with $\Rep A=\mC$ for some $G,W,u$ 
are parametrized bijectively by 7-tuples
$(G,W,H,Y,B,V,u)$ (with nondegenerate $B$). 
Upon specialization to characteristic zero, 
this is exactly the main result of \cite{EG1}. 

\begin{remark} In characteristic zero, it is known 
from \cite{De} that these are all finite dimensional
triangular Hopf algebras; see also \cite{EG2}. \end{remark}

\subsection{Taft's algebras} It is interesting to note that the principle of
proof of Theorem \ref{super} generalizes to many other situations. In this
section we consider an example when $\cC =\Rep(H_l)$ where $H_l$ is Taft's Hopf
algebra (see \cite{T}) defined as follows: choose a primitive $l-$th root of 
unity $\zeta$ (thus we assume that $\chara(k)$ does not divide $l$), then
$$H_l=\langle g,x | g^l=1, x^l=0, gxg^{-1}=\zeta x\rangle$$
$$\Delta(g)=g\otimes g, \Delta(x)=x\otimes 1+g\otimes x,$$
$$\eps(g)=1, \eps(x)=0, S(g)=g^{-1}, S(x)=-g^{-1}x.$$

\begin{example} The algebra $H_2$ is called the Sweedler's algebra. Note
that $\Rep(H_2)$ is equivalent to $\Rep(W)$ where $W$ is a one dimensional odd
vector space. \end{example}

\begin{theorem} For any divisor $d$ of $l$ there is exactly one nonsemisimple
indecomposable exact module category over $\Rep(H_l)$ with $d$ simple 
objects and exactly one one-parameter family of semisimple indecomposable 
module categories over $\Rep(H_l)$ with $d$ simple objects.
\end{theorem}

\begin{proof} Following the proof of Theorem \ref{super} we are going to
classify simple from the right algebras $A$ in the category $\Rep(H_l)$.
An algebra $A$ in the category $\Rep(H_l)$ is the same as a usual algebra
endowed with additional structures:

(i) A multiplication-preserving action of $g\in H_l$; $g^l=1$;

(ii) An action of the element $x\in H_l$ given by an operator $\p :A\to A$
satisfying $\p(ab)=\p (a)b+g(a)\p(b)$; $\p^l=0$; $\p (g(a))=
\zeta^{-1}g(\p(a))$.

Let us introduce a filtration on $A$: $A_{-1}=0; A_i=\{ a\in A | \p (a)\in 
A_{i-1}\}$. Completely analogously to the proof of Theorem \ref{super} we
have:

(a) $\p(A_i)\subset A_{i-1}$;

(b) The filtration is $g-$invariant, $g(A_i)\subset A_i$;

(c) The filtration is multiplicative, $A_iA_j\subset A_{i+j}$;

(d) The filtration is exhausting, $A_l=A$;

(e) For any $a\in A_i\setminus A_{i-1}$, $\p(a)\in A_{i-1}\setminus A_{i-2}$;

(f) The algebra $A_0$ has no nontrivial $g-$invariant ideals.

Thus the algebra $A_0$ is of the form $k[G/H]$ where $H$ is a subgroup of
$G=\langle g\rangle$ (observe that any such subgroup is cyclic). Assume that
$H$ is the (unique) subgroup of order $d$. In the case $A=A_0$ it is easy to 
check
that the category $\Mod_{\Rep(H_l)}(A)$ is a nonsemisimple module category with
$d$ simple objects. So assume that $A\ne A_0$. Let $e_s, s\in G/H$ denote the
minimal central idempotents in $A_0$. It is easy to see that $A_1$ contains a 
unique element $y$ such that $\p(y)=1$, $y=\sum_{s\in S}e_{gs}ye_s$,
$g(y)=\zeta^{-1} y$. It is also easy to see that $A$ is generated by
$A_0$ and $y$. Now an easy calculation shows that $\p(y^m)=(1+\zeta^{-1}+
\zeta^{-2}+\ldots \zeta^{1-m})y^{m-1}$ and by induction we have $y^m\in A_m
\setminus A_{m-1}$ for $m<l$. Finally, $\p(y^l)=0$ and hence $y^l=\lambda 1$
for some $\lambda \in k$. Conversely, it is not difficult to see that the 
algebra $A(d,\lambda)$ generated by $A_0$ and $y$ with the relations above is 
a simple from the right algebra in the category $\Rep(H_l)$. We leave to the 
reader to check that the category $\Mod_{\Rep(H_l)}(A(d,\lambda))$ is 
semisimple with $d$ simple objects (note that $A(d,\lambda)$ is projective as 
an object of $\Rep(H_l)$). It is obvious that the algebras $A(d,\lambda)$ are
pairwise nonisomorphic. Using the fact that all simple objects in $\Rep(H_l)$
are invertible one shows that for an indecomposable exact module category
$\cM$ over $\Rep(H_l)$ the algebra $\iHom(M,M)$ with simple $M\in \cM$ does
not depend on the choice of $M$. Thus the semisimple indecomposable exact 
module categories $\Mod_{\Rep(H_l)}(A(d,\lambda))$ over $\Rep(H_l)$ are 
pairwise nonequivalent. The Theorem is proved.
\end{proof}

\begin{remark} The algebra $A(l,\lambda)$ was studied by S.~Montgomery and
H.-J.~Schneider in \cite{MS}. It would be interesting to interpret the results
of \cite{MS} in our language. \end{remark}

\end{document}